\input amstex
\input amsppt.sty

\def\im{\operatorname{im}}
\def\rank{\operatorname{rank}}
\def\Tr{\operatorname{Tr}}

\topmatter
\title
Unified treatment of multisymplectic 3-forms in dimension 6
\endtitle
\author
Jarol\'\i m Bure\v s{} and Ji\v r\'\i{} Van\v zura
\endauthor
\address
Faculty of Mathematics and Physics, Charles University, Sokolovsk\'a{}
83, 186 00 Prague, Czech Republic
\endaddress
\email
jbures\@karlin.mff.cuni.cz
\endemail
\address
Institute of Mathematics, Academy of Sciences of the Czech Republic, \v
Zi\v zkova 22, 616 62 Brno, Czech Republic
\endaddress
\email
vanzura\@ipm.cz
\endemail
\thanks
Supported by the Grant Agency of the Academy of Sciences of the Czech
Republic, grant no\. A1019204.
\endthanks
\abstract
On a 6-dimensional real vector space $V$ there are three types of
multisymplectic 3-forms. We present in this paper a unified treatment of
these three types. Forms of each type represent a subset of $\Lambda^3
V^*$. In two cases they are open subsets, in the third one it is a
submanifold of codimension 1. We study the geometry of these subsets.
\endabstract
\endtopmatter
\NoRunningHeads

\document
\head
0. Introduction
\endhead

We shall consider a 6-dimensional real vector space $V$. Let us recall
that a multisymplectic 3-form on $V$ is a 3-form $\omega$ such that the
associated homomorphism
$$
\kappa:V\rightarrow\Lambda^2V^*,\quad\kappa v=\iota_v\omega=
\omega(v,\cdot,\cdot)
$$
is injective. We denote $\Lambda_{ms}^3V^*$ the subset of $\Lambda^3V^*$
consisting of all multisymplectic forms. It is easy to see that
$\Lambda_{ms}^3V^*$ is an open subset. The natural action of $GL(V)$ on
$\Lambda^3V^*$ preserves $\Lambda_{ms}^3V^*$. It is well known that
under this action $\Lambda_{ms}^3V^*$ decomposes into three orbits (see
e\. g\. [D], [H]). Two of them are open orbits, the third one is
a submanifold of codimension 1. As representatives of these orbits we can
take the following 3-forms. (We choose a basis $e_1,\dots,e_6$ of $V$,
and we denote $\alpha_1,\dots,\alpha_6$ the corresponding dual basis.)
\roster
\item"{(1)}"$\omega_+=\alpha_1\wedge\alpha_2\wedge\alpha_3+
           \alpha_4\wedge\alpha_5\wedge\alpha_6$,
\item"{(2)}"$\omega_-=\alpha_1\wedge\alpha_2\wedge\alpha_3+
             \alpha_1\wedge\alpha_4\wedge\alpha_5+
             \alpha_2\wedge\alpha_4\wedge\alpha_6-
             \alpha_3\wedge\alpha_5\wedge\alpha_6$,
\item"{(3)}"$\omega_0=\alpha_1\wedge\alpha_4\wedge\alpha_5+
             \alpha_2\wedge\alpha_4\wedge\alpha_6+
             \alpha_3\wedge\alpha_5\wedge\alpha_6$.
\endroster
The open set containing the form $\omega_+$ ($\omega_-$) we shall denote
$U_+$ ($U_-$), and the codimension 1 submanifold containing $\omega_0$
we shall denote $U_0$. There is also another possible characterization of
these orbits. Namely, for any 3-form $\omega$ we define
$$
\Delta^2(\omega)=\{v\in V;(\iota_v\omega)\wedge(\iota_v\omega)\}=0.
$$
In other words, the subset $\Delta^2(\omega)\subset V$ consists of all
vectors $v\in V$ such that the 2-form $\iota_v\omega$ is decomposable.
A computation shows that
$$
\gather
\Delta^2(\omega_+)=[e_1,e_2,e_3]\cup [e_4,e_5,e_6],\\
\Delta^2(\omega_2)=\{0\},\\
\Delta^2(\omega_3)=[e_1,e_2,e_3].
\endgather
$$
We find easily that
\roster
\item"{(1)}" $\omega\in U_+$ if and only if $\Delta^2(\omega)$ consists
of the union of two transversal 3-dimensional subspaces.
\item"{(2)}" $\omega\in U_-$ if and only if $\Delta^2(\omega)=\{0\}$.
\item"{(3)}" $\omega\in U_0$ if and only if $\Delta^2(\omega)$ is a
3-dimensional subspace.
\endroster

We consider now a multisymplectic 3-form $\omega$, and we choose a
nonzero 6-form $\theta$ on $V$. It is easy to see that there exists
a unique endomorphism $Q:V\rightarrow V$ such that
$$
(\iota_v\omega)\wedge\omega=\iota_{Qv}\theta.\tag{*}
$$
We shall now study the form of the endomorphism $Q$.

\head
1. The product case
\endhead
Let us assume that $\omega\in U_+$. Then $\Delta^2
(\omega)=V'_3\cup V''_3$, where $V'_3$ and $V''_3$ are transversal
3-dimensional subspaces. Our main aim in this case is to prove that
after the necessary normalization the endomorphism $Q$ is a product
structure, i\. e\. it satisfies $Q^2=I$, and its associated subspaces
are the subspaces $V'_3$ and $V''_3$.

If $v\in V'_3$, $v\ne0$ then applying $\iota_v$
to (*), we get
$$
0=(\iota_v\omega)\wedge(\iota_v\omega)=\iota_v\iota_{Qv}\theta,
$$
which shows that the vectors $v$ and $Qv$ are linearly dependent.
This means that there is a function $\lambda_1:V'_3-\{0\}\rightarrow\Bbb
R$ such that $Qv=\lambda_1(v)v$ for every $v\in V'_3-\{0\}$. It is easy to
see that the function $\lambda_1$ is constant. Namely, taking two
linearly independent vectors $v_1,v_2\in V'_3$, we get
$$
\lambda_1(v_1+v_2)v_1+\lambda_1(v_1+v_2)v_2=Q(v_1+v_2)=Q(v_1)+Q(v_2)=
\lambda_1(v_1)v_1+\lambda_1(v_2)v_2,
$$
which implies that $\lambda(v_1)=\lambda(v_2)$. Consequently, we have
$Qv=\lambda_1v$ for every $v\in V'_3$. Similarly we find that there is a
constant $\lambda_2$ such that $Qv=\lambda_2v$ for every $v\in V''_3$.
Now, we are going to prove that $\lambda_1+\lambda_2=0$. We shall need
the following lemma.
\proclaim{1.1. Lemma}
If $\omega\in U_+$, $v'\in V'_3$ and $v''\in V''_3$, then $\iota_{v'}
\iota_{v''}\omega=0$.
\endproclaim
\demo{Proof}
The lemma is obvious for the form $\omega_+$. But then it holds for
every form $\omega\in U_+$.
\enddemo
Let us take two vectors $v'\in V'_3$ and $v''\in V''_3$, $v'\ne0$,
$v''\ne0$. We have
$$
(\iota_{v'}\omega)\wedge\omega=\iota_{Qv'}\theta=\lambda_1\iota_{v'}
\theta.
$$
Applying $\iota_{v''}$ to the above equation, we get
$$
\gather
(\iota_{v''}\iota_{v'}\omega)\wedge\omega+(\iota_{v'}\omega)\wedge
(\iota_{v''}\omega)=\lambda_1\iota_{v''}\iota_{v'}\theta\\
(\iota_{v'}\omega)\wedge(\iota_{v''}\omega)=\lambda_1\iota_{v''}
\iota_{v'}\theta.
\endgather
$$
Along the same lines we get
$$
\gather
(\iota_{v''}\omega)\wedge\omega=\iota_{Qv''}\theta=\lambda_2\iota_{v''}
\theta\\
(\iota_{v'}\iota_{v''}\omega)\wedge\omega+(\iota_{v''}\omega)\wedge
(\iota_{v'}\omega)=\lambda_2\iota_{v'}\iota_{v''}\theta\\
(\iota_{v''}\omega)\wedge(\iota_{v'}\omega)=\lambda_2\iota_{v'}
\iota_{v''}\theta.
\endgather
$$
From the last two results we obtain
$$
0=(\iota_{v'}\omega)\wedge(\iota_{v''}\omega)-
(\iota_{v''}\omega)\wedge(\iota_{v'}\omega)=
\lambda_1\iota_{v''}\iota_{v'}\theta-
\lambda_2\iota_{v'}\iota_{v''}\theta=
(\lambda_1+\lambda_2)\iota_{v''}\iota_{v'}\theta,
$$
which implies $\lambda_1+\lambda_2=0$. We set now $\lambda=\lambda_1=
-\lambda_2$. Obviously $\lambda\ne0$. Otherwise we would have $\Delta^2(\omega)=V$, which is a contradiction. Further, we get $Q^2=\lambda^2I$. Now we can see that
the automorphisms
$$
S_+=\frac{1}{\lambda}Q\text{ and }S_-=-\frac{1}{\lambda}Q
\text{ satisfy }S_+^2=I\text{ and }S_-^2=I,
$$
i\. e\. they define product structures on $V$, and $S_-=-S_+$. Setting
$$
\theta_+=\lambda\theta,\quad\theta_-=-\lambda\theta,
$$
we get
$$
(\iota_v\omega)\wedge\omega=\iota_{S_+v}\theta_+,\quad
(\iota_v\omega)\wedge\omega=\iota_{S_-v}\theta_-.
$$
In the sequel we shall denote $S=S_+$ and $\theta=\theta_+$. The same
results which are valid for $S_+$ hold also for $S_-$.

\proclaim{1.2. Lemma}
If $v'\in V'_3$, $v'\ne0$, then the kernel $K(\iota_{v'}\omega)$ of
the 2-form $\iota_{v'}\omega$ equals to $[v',V''_3]$. If $v''\in V''_3$,
$v''\ne0$, then the kernel $K(\iota_{v''}\omega)$ of the 2-form
$\iota_{v''}\omega$ equals to $[v'',V'_3]$.
\endproclaim
\demo{Proof}
If $v'\in V'_3$, $v'\ne0$, then the 2-form $\iota_{v'}\omega$ is a
nonzero decomposable form. Consequently $\dim K(\iota_{v'}\omega)=4$.
Obviously $v'\in K(\iota_{v'}\omega)$, and by virtue of Lemma 1.1 also
any vector from $V''$ belongs to $K(\iota_{v'}\omega)$. This proves that
$K(\iota_{v'}\omega)=[v',V''_3]$. The second assertion follows along the
same lines.
\enddemo

\proclaim{1.3. Lemma}
For any $v\in V$ there is $\iota_{Sv}\iota_v\omega=0$.
\endproclaim
\demo{Proof}
Let us assume that $S|V'_3=I$ and $S|V''_3=-I$. Then for arbitrary $v=
v'+v''$ with $v'\in V'_3$ and $v''\in V''_3$ we have
$$
\iota_{Sv}\iota_v\omega=\iota_{S(v'+v'')}\iota_{v'+v''}\omega=
\iota_{v'-v''}\iota_{v'+v''}\omega=2\iota_{v'}\iota_{v''}\omega=0.
$$
\enddemo

\proclaim{1.4. Proposition}
There exists a unique (up to the sign) product structure $S\ne I$ on $V$ such that the form $\omega$ satisfies the relation
$$
\omega(Sv_1,v_2,v_3)=\omega(v_1,Sv_2,v_3)=\omega(v_1,v_2,Sv_3)\quad
\text{for any }v_1,v_2,v_3\in V.
$$
\endproclaim
\demo{Proof}
We shall prove first that the product structure $S$ defined above
satisfies this relation. According to the above lemma we have $\iota_v\iota_{Sv}\omega=0$ for any $v\in V$. Therefore we have
$$
0=\omega(S(v_1+v_2),v_1+v_2,v_3)=\omega(Sv_1,v_2,v_3)+\omega(Sv_2,v_1,v_3),
$$
which implies
$$
\omega(Sv_1,v_2,v_3)=\omega(v_1,Sv_2,v_3).
$$
The second equality now easily follows.
Obviously, the opposite product structure $-S$ satisfies the same
relation. It remains to prove that there is no other product structure
with the same property. Let $\tilde{S}$ be another product stucture
with the above property. Then there is a unique automorphism
$A:V\rightarrow V$ such that $\tilde{S}=SA$. We have then
$$
\gather
\omega(v_1,\tilde{S}v_2,v_3)=\omega(v_1,v_2,\tilde{S}v_3)\\
\omega(v_1,SAv_2,v_3)=\omega(v_1,v_2,SAv_3)\\
\omega(Sv_1,Av_2,v_3)=\omega(Sv_1,v_2,Av_3)\\
(\iota_{Sv_1}\omega)(Av_2,v_3)=(\iota_{Sv_1}\omega)(v_2,Av_3).
\endgather
$$ 
Because $S$ is an automorphism we get the equality
$$
(\iota_{v_1}\omega)(Av_2,v_3)=(\iota_{v_1}\omega)(v_2,Av_3).
$$
Let us take a vector $v'_1\in V'_3$. Then for any $v'_2\in V'_3$ we have
$$
0=(\iota_{v'_1}\omega)(Av'_2,v'_1)=(\iota_{v'_1}\omega)(v'_2,Av'_1).
$$
Because $v'_2$ is arbitrary, we can see that $Av'_1$ belongs to the
kernel $K(\iota_{v'}\omega)$. This means that there is $\lambda(v'_1)\in\Bbb R$
and $v''\in V''_3$ such that $Av'_1=\lambda(v'_1) v'_1+v''$. Now we can easily
see that there is $\lambda\in\Bbb R$ and a homomorphism $\varphi:V'_3\rightarrow V''_3$ such that
$$
Av'_1=\lambda v'_1+\varphi v'_1
$$
for every $v'_1\in V'_1$. Similarly we find $\mu\in\Bbb R$ and
a homomorphism $\psi:V''_3\rightarrow V'_3$ such that
$$
Av''_1=\mu v''_1+\psi v''_1
$$
for every $v''_1\in V''_3$. Taking a fixed $v'_2\in V'_3$ and arbitrary
$v''_1,v''_3\in V''_3$, we get
$$
\gather
(\iota_{v''_1}\omega)(Av'_2,v''_3)=(\iota_{v''_1}\omega)(v'_2,Av''_3)\\
(\iota_{v''_1}\omega)(\varphi v'_2,v''_3)=0,\\
(\iota_{\varphi v'_2}\omega)(v''_1,v''_3)=0.
\endgather
$$
For any $v'_1,v'_3\in V'_3$ we have by virtue of Lemma 1.1
$$
(\iota_{\varphi v'_2}\omega)(v'_1,v''_3)=0,\quad
(\iota_{\varphi v'_2}\omega)(v'_1,v'_3)=0,
$$
which together with the preceding result shows that $\iota_{\varphi v'_2}
\omega=0$. The form $\omega$ is multisymplectic and consequently
$\varphi v'_2=0$. We have thus shown that $\varphi=0$. Similarly we find
that $\psi=0$. This proves that $AV'_3\subset V'_3$, $AV''_3\subset V''_3$
and that $A|V'_3=\lambda I$, $A|V''_3=\mu I$. Because $\tilde{S}^2=I$,
we find esily that $\lambda=\pm1$ and $\mu=\pm1$. Now the proof easily follows.
\enddemo

\head
2. The complex case
\endhead

In this section we present only the relevant results. Proofs can be
found in [PV].

Let $\omega$ be a 3-form on $V$ such that $\Delta^2(\omega)=\{0\}$. This
means that for any $v\in V$, $v\ne0$ there is $(\iota_v\omega)\wedge
(\iota_v\omega)\ne0$. This implies that $\rank(\iota_v\omega)\geq4$. On
the other hand obviously $\rank(\iota_v\omega)\leq4$. Consequently, for
any $v\ne0$ $\rank(\iota_v\omega)=4$. Thus the kernel $K(\iota_v\omega)$
of the 2-form $\iota_v\omega$ has dimension 2. Moreover
$v\in K(\iota_v\omega)$. We have
$$
(\iota_v\omega)\wedge\omega=\iota_{Qv}\theta.
$$
If $v\ne0$ then $(\iota_v\omega)\wedge\omega\ne0$, and this shows that
$Q$ is an automorphism. It is also obvious that if $v\ne0$, then
the vectors $v$ and $Qv$ are linearly independent (apply $\iota_v$ to
the last equality).
\proclaim{2.1. Lemma}
For any $v\in V$ there is $\iota_{Qv}\iota_v\omega=0$, i\. e\. $Qv\in
K(\iota_v\omega)$.
\endproclaim
This lemma shows that if $v\ne0$, then $K(\iota_v\omega)=[v,Qv]$. Applying
$\iota_{Qv}$ to the equality $(\iota_v\omega)\wedge\omega=\iota_{Qv}
\theta$ and using the last lemma we obtain easily the following result.
\proclaim{2.2. Lemma}
For any $v\in V$ there is $(\iota_v\omega)\wedge(\iota_{Qv}\omega)=0$.
\endproclaim

Lemma 2.1 shows that $v\in K(\iota_{Qv}\omega)$. Because $v$ and $Qv$ are
linearly independent, we can see that
$$
K(\iota_{Qv}\omega)=[v,Qv]=K(\iota_v\omega).
$$

It can be proved that that there is $\lambda\in\Bbb R$ such that
$Q^2=-\lambda^2I$. We can now see that the automorphisms
$$
J_+=\frac{1}{\lambda}Q\text{ and }J_-=-\frac{1}{\lambda}Q
\text{ satisfy }J_+^2=-I\text{ and }J_-^2=-I,
$$
i\. e\. they define complex structures on $V$, and $J_-=-J_+$.
Setting
$$
\theta_+=\lambda\theta,\quad\theta_-=-\lambda\theta
$$
we get
$$
(\iota_v\omega)\wedge\omega=\iota_{J_+v}\theta_+,\quad
(\iota_v\omega)\wedge\omega=\iota_{J_-v}\theta_-.
$$

In the sequel we shall denote $J=J_+$ and $\theta=\theta_+$. The same
results which are valid for $J_+$ hold also for $J_-$.
\proclaim{2.3. Lemma}
There exists  a unique (up to the the sign) complex structure $J$ on $V$
such that the form $\omega$ satisfies the relation
$$
\omega(Jv_1,v_2,v_3)=\omega(v_1,Jv_2,v_3)=\omega(v_1,v_2,Jv_3)\quad
\text{for any }v_1,v_2,v_3\in V.
$$
\endproclaim

\head
3. The tangent case
\endhead
Let us assume that $\omega\in U_0$. We denote $V_0=\Delta^2(\omega)$.
If $v\in V_0$, $v\ne0$, then applying $\iota_v$ to (*), we get
$$
0=(\iota_v\omega)\wedge(\iota_v\omega)=\iota_v\iota_{Qv}\theta,
$$
which shows again that the vectors $v$ and $Qv$ are linearly dependent.
Consequently, there exists a function $\lambda:V_0-\{0\}\rightarrow
\Bbb R$ such that $Qv=\lambda(v)v$ for any $v\in V_0-\{0\}$. It is easy
to see that this function is constant. We shall need the following two
lemmas.
\proclaim{3.1. Lemma}
For any $\alpha\in V^*$ we have $(\iota_v\omega)\wedge\omega\wedge\alpha=
-\alpha(Qv)\theta$.
\endproclaim
\demo{Proof}
For a fixed $\alpha\in V^*$ there exists a unique $l_{\alpha}\in V^*$
such that
$$
(\iota_v\omega)\wedge\omega\wedge\alpha=l_{\alpha}(v)\theta.
$$
Hence we get
$$
\gather
(\iota_{Qv}\theta)\wedge\alpha=l_{\alpha}(v)\theta\\
\iota_{Qv}(\theta\wedge\alpha)-\alpha(Qv)\theta=l_{\alpha}(v)\theta\\
-\alpha(Qv)\theta=l_{\alpha}(v)\theta\\
-\alpha(Qv)=l_{\alpha}(v),
\endgather
$$
which finishes the proof.
\enddemo
\proclaim{3.2. Lemma}
Let $\alpha\in V^*$ be such that $\alpha|V_0=0$. Then we have
$(\iota_v\omega)\wedge\omega\wedge\alpha=0$.
\endproclaim
\demo{Proof}
The formula can be verified for the form $\omega_0$ by a direct
computation. But then it must be true for any 3-form $\omega\in U_0$.
\enddemo
Using these two lemmas, we get for any 1-form $\alpha$ with
$\alpha|V_0=0$
$$
0=(\iota_v\omega)\wedge\omega\wedge\alpha=-\alpha(Qv)\theta,
$$
which shows that $\alpha(Qv)=0$. We have thus proved that for any $v\in
V$ we have $Qv\in V_0$, i\. e\. $\im Q\subset V_0$. Further, for any
$v\in V$ we have $Q^2v=Q(Qv)=\lambda Qv$. This shows that the endomorphism
$Q$ satisfies the equation
$$
Q(Q-\lambda I)=0.
$$

Our next aim is to prove that the above constant $\lambda$ is zero. Let
us assume on the contrary that $\lambda\ne0$. Then there are subspaces
$R_0,R_{\lambda}\subset V$ such that
$$
V=R_0\oplus R_{\lambda},\quad Q|R_0=0,Q|R_{\lambda}=\lambda I.
$$
Obviously, both these subspaces are nontrivial. $R_0\ne0$ because $\ker
Q\subset R_0$, and $R_{\lambda}\ne0$ because $R_{\lambda}\supset V_0$. On the other hand
for any $v\in R_0$ we have
$$
\gather
(\iota_v\omega)\wedge\omega=\iota_{Qv}\theta=0\\
(\iota_v\omega)\wedge(\iota_v\omega)=0.
\endgather
$$
This shows that $v\in V_0$. Consequently, we get the inclusion $R_0
\subset V_0\subset R_{\lambda}$, which is a contradiction. We have thus
proved that $\lambda=0$ and that $Q^2=0$. Because for every $v\notin V_0$
we have $Qv\ne0$ (otherwise we would have $v\in V_0$), it is easy to see
that $\im Q=\ker Q=V_0$. The endomorphisms $Q$ satisfying $Q^2=0$ are
in differential geometry usually called tangent structures, and very often they are denoted by
$T$. But because we would have here already too many $T$'s, we have
decided to introduce the notation $F=Q$. We shall call the endomorphism
$F$ tangent structure. Let us remark that when speaking about tangent
structure, we always assume that $F^2=0$ and $\im F=\ker F$.

\proclaim{3.3. Lemma}
For any $v\in V$ we have $\iota_v\iota_{Fv}\omega=0$.
\endproclaim
\demo{Proof}
We start with the equality
$$
(\iota_v\omega)\wedge\omega=\iota_{Fv}\theta.
$$
Applying $\iota_{Fv}$ we get
$$
\gather
(\iota_{Fv}\iota_v\omega)\wedge\omega+(\iota_v\omega)\wedge
(\iota_{Fv}\omega)=0\\
-(\iota_v\iota_{Fv}\omega)\wedge\omega+(\iota_v\omega)\wedge
(\iota_{Fv}\omega)=0\\
-\iota_v(\iota_{Fv}\omega\wedge\omega)+
2(\iota_v\omega)\wedge(\iota_{Fv}\omega)=0.\\
\endgather
$$
Applying $\iota_v$ we have
$$
(\iota_v\omega)\wedge(\iota_v\iota_{Fv}\omega)=0.
$$
If the 1-form $\iota_v\iota_{Fv}\omega$ were not zero, then it would
exist a 1-form $\sigma$ such that $\iota_v\omega=\sigma\wedge
\iota_v\iota_{Fv}\omega$, and we would get
$$
(\iota_v\omega)\wedge(\iota_v\omega)=\sigma\wedge
(\iota_v\iota_{Fv}\omega)\wedge\sigma\wedge(\iota_v\iota_{Fv}\omega)=0
$$
for every $v\in V$, which is a contradiction.
\enddemo
\proclaim{3.4. Lemma}
For any three vectors $v_1,v_2,v_3\in V$ we have
$$
\omega(Fv_1,v_2,v_3)=\omega(v_1,Fv_2,v_3)=\omega(v_1,v_2,Fv_3).
$$
\endproclaim
\demo{Proof}
By virtue of Lemma 3.3 we have
$$
0=\omega(v_1+v_2,F(v_1+v_2),v_3)=\omega(v_1,Fv_2,v_3)+\omega(v_2,Fv_1,v_3),
$$
which implies
$$
\omega(Fv_1,v_2,v_3)=\omega(v_1,Fv_2,v_3).
$$
The rest of the proof is easy.
\enddemo

Let us notice that the construction of the tangent struture $F$ depends
on the choice of the 6-form $\theta$. Any other nonzero 6-form is
a nonzero real multiple $a\theta$ and the relevant construction gives
the tangent structure $(1/a)F$. In other words, the 3-form $\omega\in U_0$
determines a tangent structure up to a nonzero real multiple.

We shall
now show another possibility how to obtain these tangent structures. It
is easy to see that if $v,v'$ are two vectors from the subspace
$V_0(\omega_0)=[e_1,e_2,e_3]$, then $\iota_v\iota_{v'}\omega_0=0$.
Consequently, we have the following lemma.
\proclaim{3.5. Lemma}
Let $\omega\in U_0$. Then for any two vectors $v,v'\in
V_0=\Delta^2(\omega)$ we have $\iota_v\iota_{v'}\omega=0$.
\endproclaim

\proclaim{3.6. Lemma}
Let $R_3\subset V$ be a 3-dimensional subspace such that for any two
vectors $v,v'\in R_3$ there is $\iota_v\iota_{v'}\omega=0$. Then
$R_3=\im F$.
\endproclaim
\demo{Proof}
Let $v,v'\in R_3$. Then we have
$$
\gather
(\iota_{v'}\omega)\wedge\omega=\iota_{Fv'}\theta\\
(\iota_v\iota_{v'}\omega)\wedge\omega+(\iota_{v'}\omega)\wedge
(\iota_v\omega)=\iota_v\iota_{Fv'}\theta\\
(\iota_{v'}\omega)\wedge(\iota_v\omega)=\iota_v\iota_{Fv'}\theta.
\endgather
$$
Because the left hand side of this equality is symmetric with respect to
$v$ and  $v'$, we have
$$
\gather
\iota_v\iota_{Fv'}\theta=\iota_{v'}\iota_{Fv}\theta\\
\theta(Fv',v,\cdot,\cdot,\cdot,\cdot)=
\theta(Fv,v',\cdot,\cdot,\cdot,\cdot)\\
\theta(Fv,v',\cdot,\cdot,\cdot,\cdot)=
-\theta(v,Fv',\cdot,\cdot,\cdot,\cdot)
\endgather
$$
for any two vectors $v,v'\in R_3$.

Let us assume first that $R_3\cap\im F$ is 0-dimensional. Then, taking a
basis $v_1,v_2,v_3\in R_3$, we get a basis $v_1,v_2,v_3,Fv_1,Fv_2,Fv_3$
of $V$, and consequently we have $\theta(v_1,v_2,v_3,Fv_1,Fv_2,Fv_3)
\ne0$. We take the vectors $v_1,v_2,v_3,v_1,Fv_2,Fv_3$. Applying the
last formula, we get
$$
0\ne\omega(Fv_1,v_2,v_3,v_1,Fv_2,Fv_3)=
-\omega(v_1,Fv_2,v_3,v_1,Fv_2,Fv_3)=0,
$$
which is a contradiction.

Next, let us assume that $R_3\cap\im F$ is 1-dimensional. Obviously
$FR_3$ is 2-dimensional. Then there are two possibilities.
(1) Either $FR_3\supset R_3\cap\im F$. Then there are vectors
$v_1,v_2\in R_3$ such that $v_1,v_2,Fv_1$ is a basis of $R_3$. Then we
can find a vector $v_3$ such that $v_1,v_2,Fv_1,v_3,Fv_2,Fv_3$ is a basis
of $V$. Taking the vectors $v_1,v_2,v_1,v_3,Fv_2,Fv_3$ and applying
the above formula, we get
$$
0\ne\theta(Fv_1,v_2,v_1,v_3,Fv_2,Fv_3)=-\theta(v_1,Fv_2,v_1,v_3,Fv_2,Fv_3)=0,
$$  
which is a contradiction. (2) Or $(FR_3)\cap(R_3\cap\im F)=0$. Then we can
take a basis of $R_3$ in the form $v_1,v_2,Fv_3$, and we can complete it
to a basis $v_1,v_2,Fv_3,Fv_1,Fv_2,v_3$ of $V$. This time we take
the vectors $v_1,v_2,Fv_3,v_1,Fv_2,v_3$ and we apply the same formula.
$$
0\ne\theta(Fv_1,v_2,Fv_3,v_1,Fv_2,v_3)=-\theta(v_1,Fv_2,Fv_3,v_1,Fv_2,v_3)
=0,
$$
which is again a contradiction.

It remains to consider the case when $R_3\cap\im F$ is 2-dimensional.
Then there are again two possibilities. (1) Either $(FR_3)\cap(R_3\cap\im F)
\ne0$. Then we can take a basis of $R_3$ in the form $v_1,Fv_1,Fv_2$,
and we can complete it to a basis $v_1,Fv_1,Fv_2,v_2,v_3,Fv_3$. We take
the vectors $v_1,v_2,v_1,v_3,Fv_2,Fv_3$ and we apply again the formula.
$$
0\ne\theta(Fv_1,v_2,v_1,v_3,Fv_2,Fv_3)=-\theta(v_1,Fv_2,v_1,v_3,Fv_2,Fv_3)=0,
$$
which is a contradiction. (2) Or $(FR_3)\cap(R_3\cap\im
F)=0$. Then we take a basis of $R_3$ in the form $v_1,Fv_2,Fv_3$, and we
complete it to a basis $v_1,Fv_2,Fv_3,Fv_1,v_2,v_3$. Then, taking the
vectors $v_1,Fv_2,Fv_3,v_1,v_2,v_3$ we get in the same way as above
$$
0\ne\omega(Fv_1,Fv_2,Fv_3,v_1,v_2,v_3)=
-\omega(v_1,F^2v_2,Fv_3,v_1,v_2,v_3)=0,
$$
and we get again a contradiction. In this way we have proved that $R_3=
\im F$.
\enddemo

\proclaim{3.7. Lemma}
Let $\tilde{F}:V\rightarrow V$ be a tangent structure (i\. e\. an
endomorphism satisfying $\tilde{F}^2=0$ and $\im\tilde{F}=
\ker\tilde{F}$) such that
$$
\omega(\tilde{F}v_1,v_2,v_3)=\omega(v_1,\tilde{F}v_2,v_3)=
\omega(v_1,v_2,\tilde{F}v_3).
$$
Then $\im\tilde{F}=\im F$.
\endproclaim
\demo{Proof}
It suffices to prove that the 3-dimensional subspace $\im\tilde{F}$ has
the property described in the preceding lemma. Any two vectors $v,v'\in
\im\tilde{F}$ can be expressed in the form $v=\tilde{F}w$, $v'=
\tilde{F}w'$. Then we have
$$
\iota_v\iota_{v'}\omega=\iota_{\tilde{F}w}\iota_{\tilde{F}w'}\omega=
\omega(\tilde{F}w',\tilde{F}w,\cdot)=\omega(\tilde{F}^2w',w,\cdot)=0.
$$
\enddemo
\proclaim{3.8. Proposition}
Let $\omega\in U_0$. Then there exists (up to a nonzero multiple) a unique
tangent structure $F$ such that
$$
\omega(Fv_1,v_2,v_3)=\omega(v_1,Fv_2,v_3)=\omega(v_1,v_2,Fv_3)
$$
for all $v_1,v_2,v_3\in V$.
\endproclaim
\demo{Proof}
Let $F$ and $\tilde{F}$ be two tangent structures with the above
property. We introduce on $V$ two 3-forms by setting
$$
\sigma_F(v_1,v_2,v_3)=\omega(Fv_1,v_2,v_3),\quad
\sigma_{\tilde{F}}(v_1,v_2,v_3)=\omega(\tilde{F}v_1,v_2,v_3).
$$
Because by virtue of the preceding lemma $V_0=\ker F=\ker\tilde{F}$, it is
obvious that if $v\in V_0$, then $\iota_v\sigma_F=0$ and
$\iota_v\sigma_{\tilde{F}}=0$. This implies that there exist two 3-forms
$s_F$ and $s_{\tilde{F}}$ on $V/V_0$ such that
$$
\sigma_F=\pi^*s_F,\quad \sigma_{\tilde{F}}=\pi^*s_{\tilde{F}},
$$
where $\pi:V\rightarrow V/V_0$ is the projection. The tangent structures
$F$ and $\tilde{F}$ induce isomorphisms
$$
f:V/V_0\rightarrow V_0,\quad \tilde{f}:V/V_0\rightarrow V_0.
$$
We denote $A:V/V_0\rightarrow V/V_0$ the automorphism
$A=f^{-1}\tilde{f}$. For any three vectors $v_1,v_2,v_3\in V$ we find
$$
s_{\tilde{F}}(\pi v_1,\pi v_2,\pi v_3)=\sigma_{\tilde{F}}(v_1,v_2,v_3)=
\omega(\tilde{F}v_1,v_2,v_3)=\omega(\tilde{f}\pi v_1,v_2,v_3).
$$
We remind that the last term makes sense because $\tilde{f}\pi v_1\in
V_0$. Further we have
$$
\omega(\tilde{f}\pi v_1,v_2,v_3)=\omega(fA\pi v_1,v_2,v_3).
$$
Let us choose an element $w_1\in V$ such that $\pi w_1=A\pi v_1$. Then we
get
$$
\gather
\omega(fA\pi v_1,v_2,v_3)=\omega(f\pi w_1,v_2,v_3)=\omega(Fw_1,v_2,v_3)=\\
=\sigma_F(w_1,v_2,v_3)=s_F(A\pi v_1,\pi v_2,\pi v_3).
\endgather
$$
Proceeding in this way we obtain the relations
$$
\gather
s_{\tilde{F}}(\pi v_1,\pi v_2,\pi v_3)=
s_F(A\pi v_1,\pi v_2,\pi v_3),\\
s_{\tilde{F}}(\pi v_1,\pi v_2,\pi v_3)=
s_F(\pi v_1,A\pi v_2,\pi v_3),\\
s_{\tilde{F}}(\pi v_1,\pi v_2,\pi v_3)=s_F(\pi v_1,\pi v_2,A\pi v_3),
\endgather
$$
and the relation
$$
s_F(A\pi v_1,\pi v_2,\pi v_3)=s_F(\pi v_1,A\pi v_2,\pi v_3)=
s_F(\pi v_1,\pi v_2,A\pi v_3).
$$
Because the 3-form $s_F$ is nontrivial and because the homomorphism
$\kappa:V\rightarrow\Lambda^2V^*$ induces an isomorphism
$\kappa_0:V_0\rightarrow \Lambda^2(V/V_0)^*$, we can see that for any
2-form $\alpha$ on $V/V_0$ and any two vectors $z_1,z_2\in V/V_0$ we have
$$
\alpha(Az_1,z_2)=\alpha(z_1,Az_2).
$$
Let now $z\in V/V_0$ be arbitrary, and let us take 1-forms
$\beta_1,\beta_2\in(V/V_0)^*$ such that $\beta_1(z)=\beta_2(z)=0$. We shall
consider the 2-form $\beta_1\wedge\beta_2$. For any vector $z'\in V/V_0$ we have
$$
(\beta_1\wedge\beta_2)(Az,z')=(\beta_1\wedge\beta_2)(z,Az')=0,
$$
which shows that there is $\lambda(z)\in\Bbb R$ such that $Az=\lambda(z)z$.
Moreover, it can be easily seen that the function $\lambda(z)$ is a nonzero
constant. We thus get $A=\lambda I$ and this finishes the proof.
\enddemo

Choosing a nonzero 3-form $\eta\in\Lambda^3 (V/V_0)^*$, we can define an isomorphism $V/V_0\rightarrow\Lambda^2(V/V_0)^*$ by $w\mapsto\iota_w\eta$. Similarly, the monomorphism $\kappa:V\rightarrow\Lambda^2V^*$, $\kappa v=\iota_v\omega$ induces an isomorphism $\kappa_0:V_0\rightarrow\Lambda^2(V/V_0)^*$. We take now the following chain of
homomorphisms
$$
V\overset{\pi}\to\rightarrow V/V_0\rightarrow\Lambda^2(V/V_0)^*
\overset{\kappa_0^{-1}}\to\rightarrow V_0.
$$
We denote this composition by $C$.
\proclaim{3.9. Lemma}
The homomorphism $C$ is a tangent structure satisfying $C^2=0$,
$\im C=\ker C$ and the relation 
$$
\omega(Cv_1,v_2,v_3)=\omega(v_1,Cv_2,v_3)=\omega(v_1,v_2,Cv_3)
$$
for every $v_1,v_2,v_3\in V$.
\endproclaim
\demo{Proof}
Let us take any tangent structure $F$ with the above properties, and
let us define a 3-form $\sigma_F(v_1,v_2,v_3)=\omega(Fv_1,v_2,v_3)$
as before. There is a unique 3-form $s_F$ on $V/V_0$ such that
$\sigma_F=\pi^*s_F$, where $\pi:V\rightarrow V/V_0$ is the projection.
Obviously there is a nonzero $a\in\Bbb R$ such that $\eta=as_F$. For
any $v,v',v''\in V$ we have
$$
\gather
\omega(Cv,v',v'')=\eta(\pi v,\pi v',\pi v'')=as_F(\pi v,\pi v',\pi
v'')=\\
=a\omega(Fv,v',v'')=\omega(aFv,v',v''),
\endgather
$$
which shows that $C=aF$. This finishes the proof.
\enddemo

\head
4. Orbit of forms of the product type
\endhead

This is the orbit $U_+$, which represents an open submanifold in
$\Lambda^3V^*$. We take a point $\zeta\in U_+$. For the tangent space at
this point we have $T_{\zeta}U_+=\Lambda^3V^*$. Obviously, fixing a volume form $\theta_0$
on $V$, we can choose for each $\zeta\in U_+$ an appropriate volume form $\theta(\zeta)$
(out of the two differring by the sign) such that
$\theta(\zeta)=a\theta_0$ with $a>0$.
This means that we choose at the same time at each point $\zeta\in U_+$ a product
structure $P(\zeta)\in Aut(V)$. In other words, we can consider over
$U_+$ a trivial vector bundle $\Cal V$ with the fiber $V$, and on this
vector bundle we have a tensor field $P$ of type $(1,1)$ satisfying
$P^2=I$, $\dim\ker(P-I)=3$, and $\dim\ker(P+I)=3$. Our aim is to define
a product structure on $T_{\zeta}U_+$. We shall try to define such
a product structure by the formula
$$
\gather
(\Cal P(\zeta)\Omega)(v_1,v_2,v_3)=a\Omega(Pv_1,Pv_2,Pv_3)+\\
+b[\Omega(Pv_1,Pv_2,v_3)+\Omega(Pv_1,v_2,Pv_3)+\Omega(v_1,Pv_2,Pv_3)]+\\
+c[\Omega(Pv_1,v_2,v_3)+\Omega(v_1,Pv_2,v_3)+\Omega(v_1,v_2,Pv_3)]+\\
+d\Omega(v_1,v_2,v_3)
\endgather
$$
for any $\Omega\in T_{\zeta}U_+$. Here $P$ denotes $P(\zeta)$. It is a matter of computation 
to prove

\proclaim{4.1. Proposition}
$\Cal P(\zeta)$ satisfies $\Cal P(\zeta)^2=\Cal I$ if and only if the
quadruple $(a,b,c,d)$ is equal to one of the following 16 quadruples
$$
\gather
\big(\pm1,0,0,0),\quad(\pm\frac{1}{2},0,\mp\frac{1}{2},0\big),\quad
(0,\pm\frac{1}{2},0,\mp\frac{1}{2}),\quad(0,0,0,\pm1),\\
(\frac{1}{4},\frac{1}{4},\frac{1}{4},-\frac{3}{4}),\quad
(\frac{1}{4},-\frac{1}{4},\frac{1}{4},\frac{3}{4}),\quad
(-\frac{1}{4},-\frac{1}{4},-\frac{1}{4},\frac{3}{4}),\quad
(-\frac{1}{4},\frac{1}{4},-\frac{1}{4},-\frac{3}{4}),\\
(\frac{3}{4},\frac{1}{4},-\frac{1}{4},\frac{1}{4}),\quad
(\frac{3}{4},-\frac{1}{4},-\frac{1}{4},-\frac{1}{4}),\quad
(-\frac{3}{4},\frac{1}{4},\frac{1}{4},\frac{1}{4}),\quad
(-\frac{3}{4},-\frac{1}{4},\frac{1}{4},-\frac{1}{4}).
\endgather
$$
\endproclaim

We can define subbundles
$$
\Cal V_1=\ker (P-I),\quad\Cal V_2=\ker(P+I)
$$
satisfying $\Cal V=\Cal V_1\oplus\Cal V_2$. This decomposition enables to
introduce in the standard way forms of type $(r,s)$. We denote by the
symbol $\Cal D^{r,s}$ the subbundle of the bundle $\Lambda^*\Cal V$
consisting of forms of type $(r,s)$. Now, it is obvious that the tangent
bundle of $U_+$ can be expressed as a direct sum of four subbundles
(distributions)
$$
TU_+=\Cal D^{3,0}\oplus\Cal D^{2,1}\oplus\Cal D^{1,2}\oplus\Cal D^{0,3},
$$
where $\dim\Cal D^{3,0}=\dim\Cal D^{0,3}=1$, $\dim\Cal D^{2,1}=
\dim\Cal D^{1,2}=9$. Let us denote $\pi_1:\Cal V\rightarrow\Cal V_1$ and
$\pi_2:\Cal V\rightarrow\Cal V_2$ the projections. If $\zeta\in U_+$, we
can define vectors $\zeta_1,\zeta_2\in T_{\zeta}U_+$ by the formulas
$$
\zeta_1=\pi_1^*(\zeta|\Cal V_{1\zeta}),\quad
\zeta_2=\pi_2^*(\zeta|\Cal V_{2\zeta}).
$$
Now we can define vector fields $\omega$, $\omega_1$ and $\omega_2$ on $U_+$ by
$\omega_{\zeta}=\zeta$, $\omega_{1\zeta}=\zeta_1$ and $\omega_{2\zeta}=\zeta_2$. 
Obviously, $\omega=\omega_1+\omega_2$.

To each quadruple $(a,b,c,d)$ from Proposition 4.1 there correspond a
product structure $\Cal P$ and a subbundle $\Cal V_1=\ker(\Cal P-\Cal
I)$. Routine considerations show that the correspondence $(a,b,c,d)
\mapsto\Cal V_1$ is the following one.
$$
\alignat2
&(1,0,0,0)\mapsto\Cal D^{3,0}\oplus\Cal D^{1,2}\qquad&
&(-1,0,0,0)\mapsto\Cal D^{2,1}\oplus\Cal D^{0,3}\\
&(\frac{1}{2},0,-\frac{1}{2},0)\mapsto\Cal D^{1,2}\oplus\Cal
D^{0,3}\qquad&
&(-\frac{1}{2},0,\frac{1}{2},0)\mapsto\Cal D^{3,0}\oplus\Cal D^{2,1}\\
&(0,\frac{1}{2},0,-\frac{1}{2})\mapsto\Cal D^{3,0}\oplus\Cal
D^{0,3}\qquad&
&(0,-\frac{1}{2},0,\frac{1}{2})\mapsto\Cal D^{2,1}\oplus\Cal D^{1,2}\\
&(0,0,0,1)\mapsto\Cal D^{3,0}\oplus\Cal D^{2,1}\oplus\Cal D^{1,2}\oplus
\Cal D^{0,3}\qquad&
&(0,0,0,-1)\mapsto 0\\
&(\frac{1}{4},\frac{1}{4},\frac{1}{4},-\frac{3}{4})\mapsto\Cal D^{3,0}
\qquad&
&(\frac{1}{4},-\frac{1}{4},\frac{1}{4},\frac{3}{4})\mapsto\Cal D^{3,0}
\oplus\Cal D^{2,1}\oplus\Cal D^{1,2}\\
&(-\frac{1}{4},-\frac{1}{4},-\frac{1}{4},\frac{3}{4})\mapsto\Cal D^{2,1}
\oplus\Cal D^{1,2}\oplus\Cal D^{0,3}\qquad&
&(-\frac{1}{4},\frac{1}{4},-\frac{1}{4},-\frac{3}{4})\mapsto\Cal
D^{0,3}\\
&(\frac{3}{4},\frac{1}{4},-\frac{1}{4},\frac{1}{4})\mapsto\Cal
D^{3,0}\oplus\Cal D^{1,2}\oplus\Cal D^{0,3}\qquad&
&(\frac{3}{4},-\frac{1}{4},-\frac{1}{4},-\frac{1}{4})\mapsto\Cal
D^{1,2}\\
&(-\frac{3}{4},\frac{1}{4},\frac{1}{4},\frac{1}{4})\mapsto\Cal
D^{3,0}\oplus\Cal D^{2,1}\oplus\Cal D^{0,3}\qquad&
&(-\frac{3}{4},-\frac{1}{4},\frac{1}{4},-\frac{1}{4})\mapsto\Cal D^{2,1}
\endalignat
$$

In the sequel we are going to investigate the integrability of all these 
distributions. Our first result is easy because the distributions 
$\Cal D^{3,0}$ and $\Cal D^{0,3}$ are 1-dimensional.

\proclaim{4.2. Proposition}
The distribution $\Cal D^{3,0}$ ($\Cal D^{3,0}$) is generated by the vector field 
$\omega_1$ ($\omega_2$). The distributions $\Cal D^{3,0}$ and $\Cal D^{0,3}$ are 
integrable.
\endproclaim

Now we shall introduce on $U_+$ a flat connection $\nabla$, which is the
restriction of the canonical connection on the vector space
$\Lambda^3V^*$. Notice that for any vector field $\Omega$ on $U_+$ we have 
$\nabla_{\Omega}\omega=\Omega$. We shall need the following three lemmas.

\proclaim{4.3. Lemma}
Let $\tilde{\Omega}$ be a vector field on $U_+$ belonging to $\Cal D^{3,0}$
($\Cal D^{2,1}$, $\Cal D^{1,2}$, $\Cal D^{0,3}$). Further, let $\Omega$
be arbitrary vector field on $U_+$. Then
$$
\nabla_{\Omega}\tilde{\Omega}\in\Cal D^{3,0}\oplus\Cal D^{2,1}\quad
(\Cal D^{3,0}\oplus\Cal D^{2,1}\oplus\Cal D^{1,2},
\Cal D^{2,1}\oplus\Cal D^{1,2}\oplus\Cal D^{0,3},
\Cal D^{1,2}\oplus\Cal D^{0,3}).
$$
\endproclaim
\demo{Proof}
Let $\Theta$ be a section of the trivial vector bundle $\Cal V^*$ over
$U_+$. Then for any vector field $\Omega$ on $U_+$ we have
$$
\nabla_{\Omega}\Theta\in\Cal D^{1,0}\oplus\Cal D^{0,1}.
$$
Now the assertion of the lemma easily follows.
\enddemo

\proclaim{4.4. Lemma}
If $\Omega$ belongs to the distribution $\Cal D^{3,0}$ ($\Cal D^{0,3}$),
then we have
$$
\nabla_{\Omega}\omega_1=\Omega,\quad\nabla_{\Omega}\omega_2=0
\quad\quad
(\nabla_{\Omega}\omega_1=0,\quad\nabla_{\Omega}\omega_2=\Omega).
$$
If $\Omega$ belongs to the distribution $\Cal D^{2,1}$ ($\Cal D^{1,2}$),
then we have again
$$
\nabla_{\Omega}\omega_1=\Omega,\quad\nabla_{\Omega}\omega_2=0
\quad\quad
(\nabla_{\Omega}\omega_1=0,\quad\nabla_{\Omega}\omega_2=\Omega).
$$
\endproclaim
\demo{Proof}
We start with the equality $\omega_1+\omega_2=\omega$. If $\Omega$
belongs to $\Cal D^{3,0}$, then applying $\nabla_{\Omega}$ to this
equality we get
$$
\gather
\nabla_{\Omega}\omega_1+\nabla_{\Omega}\omega_2=\Omega\\
(\nabla_{\Omega}\omega_1)^{3,0}+(\nabla_{\Omega}\omega_1)^{2,1}+
(\nabla_{\Omega}\omega_2)^{1,2}+(\nabla_{\Omega}\omega_2)^{0,3}=\Omega,
\endgather
$$
where the superscripts denote the corresponding component. Because
$\Omega$ belongs to $\Cal D^{3,0}$ we obtain the first assertion. The
remaining assertions follow along the same lines.
\enddemo

\proclaim{4.5. Lemma}
A vector field $\Omega$ belongs to the distribution $\Cal D^{3,0}\oplus
\Cal D^{2,1}$ ($\Cal D^{1,2}\oplus\Cal D^{0,3}$) if and only if
$$
\nabla_{\Omega}\omega_2=0\quad(\nabla_{\Omega}\omega_1=0).
$$
\endproclaim
\demo{Proof}
If $\Omega$ belongs to $\Cal D^{3,0}\oplus \Cal D^{2,1}$ we know that
the above condition is satisfied. Conversely, let us assume that the
condition is satisfied. We have
$$
\Omega=\Omega^{3,0}+\Omega^{2,1}+\Omega^{1,2}+\Omega^{0,3},
$$
and we get
$$
\gather
0=\nabla_{\Omega}\omega_2=\nabla_{\Omega^{3,0}}\omega_2+
\nabla_{\Omega^{2,1}}\omega_2+\nabla_{\Omega^{1,2}}\omega_2+
\nabla_{\Omega^{0,3}}\omega_2=\\
=\nabla_{\Omega^{1,2}}\omega_2+\nabla_{\Omega^{0,3}}\omega_2=
\Omega^{1,2}+\Omega^{0,3},
\endgather
$$
which finishes the proof.
\enddemo

\proclaim{4.6. Proposition}
The distributions $\Cal D^{3,0}\oplus\Cal D^{2,1}$ and
$\Cal D^{1,2}\oplus\Cal D^{0,3}$ are integrable.
\endproclaim
\demo{Proof}
Let two vector fields $\Omega,\tilde{\Omega}$ belong to the distribution
$\Cal D^{3,0}\oplus\Cal D^{2,1}$. Then we have $\nabla_{\Omega}\omega_2=
\nabla_{\tilde{\Omega}}\omega_2=0$, and we obtain
$$
\nabla_{[\Omega,\tilde{\Omega}]}\omega_2=\nabla_{\Omega}
\nabla_{\tilde{\Omega}}\omega_2-\nabla_{\tilde{\Omega}}
\nabla_{\Omega}\omega_2=0
$$
because the connection $\nabla$ is flat. Along the same lines we can
prove the integrability of the distribution $\Cal D^{1,2}\oplus
\Cal D^{0,3}$.
\enddemo
The following lemma is obvious.
\proclaim{4.7. Lemma}
A vector field $\Omega$ belongs to the distribution $\Cal D^{2,1}\oplus
\Cal D^{1,2}$ if and only if $\Omega\wedge\omega=0$.
\endproclaim

\proclaim{4.8. Proposition}
The distribution $\Cal D^{2,1}\oplus\Cal D^{1,2}$ is not integrable.
\endproclaim
\demo{Proof}
Let $\Omega$ and $\tilde{\Omega}$ lie in $\Cal D^{2,1}$ and
$\Cal D^{1,2}$, respectively. Then we have $\Omega\wedge\omega=0$ and
$\tilde{\Omega}\wedge\omega=0$. Hence we obtain
$$
(\nabla_{\Omega}\tilde{\Omega})\wedge\omega+
\tilde{\Omega}\wedge\Omega=0,\quad
(\nabla_{\tilde{\Omega}}\Omega)\wedge\omega+\Omega\wedge\tilde{\Omega}=0.
$$
Substracting these two equalities, we have
$$
[\Omega,\tilde{\Omega}]\wedge\omega=2\Omega\wedge\tilde{\Omega}.
$$
Now it suffices to choose $\Omega$ and $\tilde{\Omega}$ in such a way
that $\Omega_{\zeta}\wedge\tilde{\Omega}_{\zeta}\ne0$ at some point
$\zeta\in U_+$. Then it is obvious that the bracket $[\Omega,
\tilde{\Omega}]$ does not lie in $\Cal D^{2,1}\oplus\Cal D^{1,2}$.
\enddemo

\proclaim{4.9. Proposition}
The distributions $\Cal D^{2,1}$ and $\Cal D^{1,2}$ are integrable.
\endproclaim
\demo{Proof}
Let $\Omega$ and $\tilde{\Omega}$ be two vector fields lying in $\Cal
D^{2,1}$. Proceeding in the same way as in the proof of preceding lemma
we find again
$$
[\Omega,\tilde{\Omega}]\wedge\omega=2\Omega\wedge\tilde{\Omega}.
$$
But this time $\Omega\wedge\tilde{\Omega}=0$, which shows that $[\Omega,
\tilde{\Omega}]$ lies in $\Cal D^{2,1}\oplus\Cal D^{1,2}$. Moreover, we
have
$$
\nabla_{[\Omega,\tilde{\Omega}]}\omega_2=
\nabla_{\Omega}\nabla_{\tilde{\Omega}}\omega_2-
\nabla_{\tilde{\Omega}}\nabla_{\Omega}\omega_2=0,
$$
which shows that $[\Omega,\tilde{\Omega}]$ lies in $\Cal D^{2,1}$.
\enddemo

\proclaim{4.10. Proposition}
There is $[\omega_1,\omega_2]=0$ and the distribution $\Cal D^{3,0}
\oplus\Cal D^{0,3}$ is integrable.
\endproclaim
\demo{Proof}
We have
$$
\nabla_{[\omega_1,\omega_2]}\omega_1=
\nabla_{\omega_1}\nabla_{\omega_2}\omega_1-
\nabla_{\omega_2}\nabla_{\omega_1}\omega_1=0-\nabla_{\omega_2}\omega_1=0,
$$
which shows that $[\omega_1,\omega_2]$ lies in $\Cal D^{1,2}\oplus
\Cal D^{0,3}$. Along the same lines we can show that $[\omega_1,\omega_2]$
lies in $\Cal D^{3,0}\oplus\Cal D^{2,1}$. This implies that
$[\omega_1,\omega_2]=0$ and that the distribution $\Cal D^{3,0}
\oplus\Cal D^{0,3}$ is integrable.
\enddemo

\proclaim{4.11. Proposition}
For any vector field $\Omega$ lying in $\Cal D^{1,2}$ ($\Cal D^{2,1}$)
the vector field $[\omega_1,\Omega]$ ($[\omega_2,\Omega]$) lies again in
$\Cal D^{1,2}$ ($\Cal D^{2,1}$). Consequently the distributions $\Cal
D^{3,0}\oplus\Cal D^{1,2}$ and $\Cal D^{2,1}\oplus\Cal D^{0,3}$ are
integrable.
\endproclaim
\demo{Proof}
Let us assume that $\Omega$ lies in $\Cal D^{1,2}$. Then we have
$$
\nabla_{[\omega_1,\Omega]}\omega_1=
\nabla_{\omega_1}\nabla_{\Omega}\omega_1-
\nabla_{\Omega}\nabla_{\omega_1}\omega_1=0-\nabla_{\Omega}\omega_1=0,
$$
which proves that $[\omega_1,\Omega]$ lies in $\Cal D^{1,2}\oplus\Cal
D^{0,3}$. Because $\Omega$ lies in $\Cal D^{1,2}$, there is $\Omega
\wedge\omega=0$. Applying $\nabla_{\omega_1}$ to this equality we find
that
$$
0=(\nabla_{\omega_1}\Omega)\wedge\omega+
\Omega\wedge\nabla_{\omega_1}\omega=
(\nabla_{\omega_1}\Omega)\wedge\omega+
\Omega\wedge\omega_1.
$$
Obviously $\Omega\wedge\omega_1=0$, and this shows that
$\nabla_{\omega_1}\Omega$ lies in $\Cal D^{2,1}\oplus\Cal D^{1,2}$.
But we can immediately see that
$$
[\omega_1,\Omega]=\nabla_{\omega_1}\Omega-\nabla_{\Omega}\omega_1=
\nabla_{\omega_1}\Omega.
$$
Consequently $[\omega_1,\Omega]$ lies not only in $\Cal D^{1,2}\oplus\Cal
D^{0,3}$, but also in $\Cal D^{2,1}\oplus\Cal D^{1,2}$. This implies
that $[\omega_1,\Omega]$ lies in $\Cal D^{1,2}$.
\enddemo

\proclaim{4.12. Proposition}
The distributions $\Cal D^{3,0}\oplus\Cal D^{2,1}\oplus\Cal D^{0,3}$ and
$\Cal D^{3,0}\oplus\Cal D^{1,2}\oplus\Cal D^{0,3}$ are integrable. The
distributions $\Cal D^{3,0}\oplus\Cal D^{2,1}\oplus\Cal D^{1,2}$ and
$\Cal D^{2,1}\oplus\Cal D^{1,2}\oplus\Cal D^{0,3}$ are not integrable.
\endproclaim
\demo{Proof}
The first assertion is easy to prove. Therefore, let us consider the
distribution $\Cal D^{3,0}\oplus\Cal D^{2,1}\oplus\Cal D^{1,2}$. We
shall take the same vector fields $\Omega$ lying in $\Cal D^{2,1}$ and
$\tilde{\Omega}$ lying in $\Cal D^{1,2}$ as in the proof of Proposition 4.8.
Then we have
$$
\gather
[\Omega,\tilde{\Omega}]\wedge\omega_1=(\nabla_{\Omega}\tilde{\Omega})
\wedge\omega_1-(\nabla_{\tilde{\Omega}}\Omega)\wedge\omega_1=\\
=\nabla_{\Omega}(\tilde{\Omega}\wedge\omega_1)-
\tilde{\Omega}\wedge(\nabla_{\Omega}\omega_1)-
\nabla_{\tilde{\Omega}}(\Omega\wedge\omega_1)+
\Omega\wedge(\nabla_{\tilde{\Omega}}\omega_1)=\\
=-\tilde{\Omega}\wedge\Omega=\Omega\wedge\tilde{\Omega}.
\endgather
$$
At the same point $\zeta\in U_+$ as in the proof of Proposition 4.8 we
have $\Omega_{\zeta}\wedge\tilde{\Omega}_{\zeta}\ne0$, which shows that
$[\Omega,\tilde{\Omega}]^{0,3}_{\zeta}\ne0$. This proves that the
distribution under consideration is not integrable.
\enddemo

We can summarize our results.
\proclaim{4.13. Proposition}
The distributions
$$
\gather
\Cal D^{3,0},\quad\Cal D^{2,1},\quad\Cal D^{1,2},\quad\Cal D^{0,3}\\
\Cal D^{3,0}\oplus\Cal D^{2,1},\quad\Cal D^{3,0}\oplus\Cal D^{1,2},\quad
\Cal D^{3,0}\oplus\Cal D^{0,3},\quad\Cal D^{2,1}\oplus\Cal D^{0,3},\quad
\Cal D^{1,2}\oplus\Cal D^{0,3}\\
\Cal D^{3,0}\oplus\Cal D^{2,1}\oplus\Cal D^{0,3},\quad
\Cal D^{3,0}\oplus\Cal D^{1,2}\oplus\Cal D^{0,3}
\endgather
$$
are integrable. The distributions
$$
\Cal D^{2,1}\oplus\Cal D^{1,2},\quad
\Cal D^{3,0}\oplus\Cal D^{2,1}\oplus\Cal D^{1,2},\quad
\Cal D^{2,1}\oplus\Cal D^{1,2}\oplus\Cal D^{0,3}
$$
are not integrable.
\endproclaim

\definition{4.14. Remark}
Requiring $\dim\ker(\Cal P-\Cal I)=\dim\ker(\Cal P+\Cal I)=10$ we have
only four possibilities how to define a product structure $\Cal P$. It
is easy to see that these product structures correspond to the
quadruples
$$
(1,0,0,0),\quad(-1,0,0,0),\quad(\frac{1}{2},0,-\frac{1}{2},0),\quad
(-\frac{1}{2},0,\frac{1}{2},0).
$$
Because all the distributions associated with these projectors are
integrable, in all these cases the Nijenhuis tensor $[\Cal P,\Cal P]=0$.
\enddefinition

\head
5. Orbit of forms of the complex type
\endhead

Here we shall study the orbit $U_-$, which also represents an open
submanifold in $\Lambda^3V^*$. Taking a point $\zeta\in U_-$, we have
$T_{\zeta}U_-=\Lambda^3V^*$. Fixing again a volume form $\theta_0$
on $V$, we can choose for each $\zeta\in U_-$ an appropriate volume form
$\theta(\zeta)$ (out of the two differring by the sign) such that
$\theta(\zeta)=a\theta_0$ with $a>0$. This enables us to choose at each point
$\zeta\in U_-$ a complex structure $J(\zeta)\in Aut(V)$. In other words,
this time we have on the trivial vector bundle $\Cal V$ a tensor field
$J$ of type $(1,1)$ satisfying $J^2=-I$. We shall again try to define a
complex structure on $T_{\zeta}U_-$ by the formula
$$
\gather
(\Cal J(\zeta)\Omega)(v_1,v_2,v_3)=a\Omega(Jv_1,Jv_2,Jv_3)+\\
+b[\Omega(Jv_1,Jv_2,v_3)+\Omega(Jv_1,v_2,Jv_3)+\Omega(v_1,Jv_2,Jv_3)]+\\
+c[\Omega(Jv_1,v_2,v_3)+\Omega(v_1,Jv_2,v_3)+\Omega(v_1,v_2,Jv_3)]+\\
+d\Omega(v_1,v_2,v_3)
\endgather
$$
for any $\Omega\in T_{\zeta}U_-$.

\proclaim{5.1. Proposition}
$\Cal J(\zeta)$ satisfies $\Cal J(\zeta)^2=-\Cal I$ if and only if the
quadruple $(a,b,c,d)$ is equal to one of the following 4 quadruples
$$
(\pm1,0,0,0),\quad(\pm\frac{1}{2},0,\pm\frac{1}{2},0).
$$
\endproclaim
The proof is a simple computation and will be omitted. We shall denote
$$
\gather
(\Cal J_1(\zeta)\Omega)(v_1,v_2,v_3)=\Omega(J(\zeta)v_1,J(\zeta)v_2,
                                            J(\zeta)v_3)\\
(\Cal J_2(\zeta)\Omega)(v_1,v_2,v_3)=\frac{1}{2}\Omega(J(\zeta)v_1,
J(\zeta)v_2,J(\zeta)v_3)+\\
+\frac{1}{2}[\Omega(J(\zeta)v_1,v_2,v_3)+
\Omega(v_1,J(\zeta)v_2,v_3)+\Omega(v_1,v_2,J(\zeta)v_3)].
\endgather
$$
The mapping $\zeta\in U_-\mapsto J_1(\zeta)$ (resp\. $\zeta\in U_-
\mapsto J_2(\zeta)$) defines an almost complex structure $\Cal J_1$
(resp\. $\Cal J_2$) on the orbit $U_-$.

\proclaim{5.2. Proposition}
The almost complex structure $\Cal J_2$ is integrable.
\endproclaim
\demo{Proof}
We denote again by $\nabla$ the canonical connection on $\Lambda^3V^*$.
Let $\Omega$ and $\tilde{\Omega}$ be two vector fields on
$U_-$. Applying $\nabla_{\tilde{\Omega}}$ to the identity $J^2=-I$, we
get
$$
(\nabla_{\tilde{\Omega}}J)J+J(\nabla_{\tilde{\Omega}}J)=0.
$$
Further, we shall use the identity
$$
\omega(Jv_1.v_2,v_3)=\omega(v_1,Jv_2,v_3),
$$
and apply to it the covariant derivative $\nabla_{\tilde{\Omega}}$. We
obtain
$$
\tilde{\Omega}(Jv_1,v_2,v_3)+\omega((\nabla_{\tilde{\Omega}}J)v_1,
v_2,v_3)=
\tilde{\Omega}(v_1,Jv_2,v_3)+\omega(v_1,(\nabla_{\tilde{\Omega}}J)v_2,
v_3).
$$
Substituing now $Jv_2$ instead of $v_2$ and $(\nabla_{\Omega}J)v_3$ instead
of $v_3$, we get the relation
$$
\gather
\tilde{\Omega}(Jv_1,Jv_2,(\nabla_{\Omega}J)v_3)=\\
-\tilde{\Omega}(v_1,v_2,(\nabla_{\Omega}J)v_3)
-\omega((\nabla_{\tilde{\Omega}}J)v_1,Jv_2,(\nabla_{\Omega}J)v_3)
-\omega(Jv_1,(\nabla_{\tilde{\Omega}}J)v_2,(\nabla_{\Omega}J)v_3).
\endgather
$$
Similarly we obtain the relations
$$
\gather
\tilde{\Omega}(Jv_1,(\nabla_{\Omega}J)v_2,Jv_3)=\\
-\tilde{\Omega}(v_1,(\nabla_{\Omega}J)v_2,v_3)
-\omega(Jv_1,(\nabla_{\Omega}J)v_2,(\nabla_{\tilde{\Omega}}J)v_3)
-\omega((\nabla_{\tilde{\Omega}}J)v_1,(\nabla_{\Omega}J)v_2,Jv_3),\\
\tilde{\Omega}((\nabla_{\Omega}J)v_1,Jv_2,Jv_3)=\\
-\tilde{\Omega}((\nabla_{\Omega}J)v_1,v_2,v_3)
-\omega((\nabla_{\Omega}J)v_1,(\nabla_{\tilde{\Omega}}J)v_2,Jv_3)
-\omega((\nabla_{\Omega}J)v_1,Jv_2,(\nabla_{\tilde{\Omega}}J)v_3).
\endgather
$$
Let us compute now
$$
\gather
2(\nabla_{\Omega}(\Cal J\tilde{\Omega}))(v_1,v_2,v_3)
=2\nabla_{\Omega}((\Cal J\tilde{\Omega})(v_1,v_2,v_3))
=\nabla_{\Omega}(\tilde{\Omega}(Jv_1,Jv_2,Jv_3)+\\
+[\tilde{\Omega}(Jv_1,v_2,v_3)+\tilde{\Omega}(v_1,Jv_2,v_3)+
  \tilde{\Omega}(v_1,v_2,Jv_3)])=
(\nabla_{\Omega}\tilde{\Omega})(Jv_1,Jv_2,Jv_3)+\\
+(\nabla_{\Omega}\tilde{\Omega})(Jv_1,v_2,v_3)+
(\nabla_{\Omega}\tilde{\Omega})(v_1,Jv_2,v_3)+
(\nabla_{\Omega}\tilde{\Omega})(v_1,v_2,Jv_3)+\\
+\tilde{\Omega}((\nabla_{\Omega}J)v_1,Jv_2,Jv_3)+
\tilde{\Omega}(Jv_1,(\nabla_{\Omega}J)v_2,Jv_3)+
\tilde{\Omega}(Jv_1,Jv_2,(\nabla_{\Omega}J)v_3)+\\
+\tilde{\Omega}((\nabla_{\Omega}J)v_1,v_2,v_3)+
\tilde{\Omega}(v_1,(\nabla_{\Omega}J)v_2,v_3)+
\tilde{\Omega}(v_1,v_2,(\nabla_{\Omega}J)v_3)=\\
=2(\Cal J\nabla_{\Omega}\tilde{\Omega})(v_1,v_2,v_3)-\\
-\omega((\nabla_{\Omega}J)v_1,(\nabla_{\tilde{\Omega}}J)v_2,Jv_3)
-\omega((\nabla_{\Omega}J)v_1,Jv_2,(\nabla_{\tilde{\Omega}}J)v_3)\\
-\omega(Jv_1,(\nabla_{\Omega}J)v_2,(\nabla_{\tilde{\Omega}}J)v_3)
-\omega((\nabla_{\tilde{\Omega}}J)v_1,(\nabla_{\Omega}J)v_2,Jv_3)\\
-\omega((\nabla_{\tilde{\Omega}}J)v_1,Jv_2,(\nabla_{\Omega}J)v_3)
-\omega(Jv_1,(\nabla_{\tilde{\Omega}}J)v_2,(\nabla_{\Omega}J)v_3).
\endgather
$$
Here we have used the previous relations. Let us notice that the expression
consisting of the last six terms is symmetric with respect to $\Omega$
and $\tilde{\Omega}$. Consequently we obtain
$$
\nabla_{\Omega}(\Cal J\tilde{\Omega})-
\nabla_{\tilde{\Omega}}(\Cal J\Omega)=
\Cal J(\nabla_{\Omega}\tilde{\Omega}-\nabla_{\tilde{\Omega}}\Omega)=
\Cal J[\Omega,\tilde{\Omega}].
$$
Writing $\Cal J\Omega$ instead of $\Omega$, we get
$$
\nabla_{\Cal J\Omega}(\Cal J\tilde{\Omega})=
-\nabla_{\tilde{\Omega}}\Omega+\Cal J[\Cal J\Omega,\tilde{\Omega}].
$$
Interchanging $\Omega$ and $\tilde{\Omega}$ we get the relation
$$
\nabla_{\Cal J\tilde{\Omega}}(\Cal J\Omega)=
-\nabla_{\Omega}\tilde{\Omega}+\Cal J[\Cal J\tilde{\Omega},\Omega].
$$
Substracting these last two relations we obtain
$$
\gather
[\Cal J\Omega,\Cal J\tilde{\Omega}]=[\Omega,\tilde{\Omega}]
+\Cal J[\Cal J\Omega,\tilde{\Omega}]
-\Cal J[\Cal J\tilde{\Omega},\Omega]\\
[\Cal J\Omega,\Cal J\tilde{\Omega}]-[\Omega,\tilde{\Omega}]
-\Cal J[\Cal J\Omega,\tilde{\Omega}]
-\Cal J[\Omega,\Cal J\tilde{\Omega}]=0,
\endgather
$$
which shows that the Nijenhuis tensor $[\Cal J,\Cal J]=0$.
\enddemo
\definition{5.3. Remark}
The almost complex structure $\Cal J_2$ was introduced in quite
different way by N.~Hitchin in [H]. He also proved the integrability and
some other properties of $\Cal J_2$.
\enddefinition

\head
6. Orbit of forms of the tangent type
\endhead

Here we shall investigatethe the last orbit $U_0$, which represents a
submanifold of codimension 1 in $\Lambda^3V^*$.
Let $\zeta\in U_0$ be arbitrary point, and let us denote $V_0(\zeta)=
\Delta^2(\zeta)$. We shall introduce three subspaces
$\Cal D_i(\zeta)\subset V$, $i=1,2,3$ in the following way:
$$
\Cal D_i(\zeta)=\{\Omega\in T_{\zeta}U_0;\Omega(v_1,v_2,v_3)=0\text
{ if the vectors }v_1,\dots,v_i\text{ belong to }V_0(\zeta)\}.
$$
It is easy to verify that $\dim\Cal D_1=1$, $\dim\Cal D_2=10$, $\dim
\Cal D_3=19$. Moreover, it is obvious that
$$
\Cal D_1\subset\Cal D_2\subset\Cal D_3.
$$

We describe first the tangent spaces to the orbit $U_0$. It is obvious
that the projection
$$
\pi_{\zeta}:GL(6,\Bbb R)\rightarrow U_0,\quad\pi_{\zeta}(\varphi)=
\varphi^*\zeta
$$
admits a smooth local section $\sigma$ defined on an open neighborhood
$W$ of $\zeta$ and such that $\sigma(\zeta)=1$. For any $\omega\in W$
we have then
$$
\omega=\sigma(\omega)^*\zeta.
$$
Let $\gamma:(-\varepsilon,\varepsilon)\rightarrow W$ be a smooth curve
such that $\gamma(0)=\zeta$. We have then
$$
\gather
\gamma(t)=\sigma(\gamma(t))^*\zeta\\
\gamma(t)(v_1,v_2,v_3)=\zeta(\sigma(\gamma(t))v_1,\sigma(\gamma(t))v_2,
\sigma(\gamma(t))v_3),
\endgather
$$
where $v_1,v_2,v_3\in V$ are arbitrary. Differentiating the last
equality at $t=0$, we get
$$
\Omega(v_1,v_2,v_3)=\zeta(Av_1,v_2,v_3)+\zeta(v_1,Av_2,v_3)+
\zeta(v_1,v_2,Av_3),
$$
where $\Omega=(d/dt)_{t=0}\gamma(t)$ and $A=(d/dt)_{t=0}\sigma(\gamma(t))$.

\proclaim{6.1. Proposition}
There is $T_{\zeta}U_0=\Cal D_3(\zeta)$.
\endproclaim
\demo{Proof}
If $\Omega\in T_{\zeta}U_0$, then according to the above formula there
is $\Omega\in\Cal D_3(\zeta)$ because $\zeta(v,v',v'')=0$ if two entries
belong to $V_0(\zeta)$. We have therefore $T_{\zeta}U_0\subset\Cal
D_3(\zeta)$. Because $\dim T_{\zeta}U_0=19$ and $\dim\Cal D_3(\zeta)=19$,
we get $T_{\zeta}U_0=\Cal D_3(\zeta)$.
\enddemo

It is obvious that it makes no sense to use in the future the notation
$\Cal D_3(\zeta)$. The following lemma can be easily verified for the
form $\omega_0$. But then it necessarily holds for any form $\zeta\in
U_0$
\proclaim{6.2. Lemma}
There is
$$
\gather
\Cal D_2(\zeta)=\{\Omega\in T_{\zeta}U_0;\Omega\wedge(\iota_v\zeta)=0
\text{ for every }v\in V_0(\zeta)\}=\\
=\{\Omega\in T_{\zeta}U_0;\Omega\wedge\beta\wedge\beta'=0
\text{ for any }\beta,\beta\in V^*\text {such that }\beta|V_0(\zeta)=
\beta'|V_0(\zeta)=0\}.
\endgather
$$
\endproclaim

On $U_0$ we have the trivial 6-dimensional vector bundle $\Cal V$ with
fiber $V$, and we can define a 3-dimensional vector subbundle $\Cal V_0$
whose fiber at $\zeta$ is $V_0(\zeta)$. We denote $\Cal W$ the
3-dimensional quotient vector bundle $\Cal V/\Cal V_0$. Moreover,
assigning to each point $\zeta\in U_0$ the vector space $\Cal
D_i(\zeta)$, we obtain over $U_0$ a vector bundle $\Cal D_i$, $i=1,2$.
In other words we have two distributions $\Cal D_1\subset\Cal D_2\subset
TU_0$. Furthermore, we have on $U_0$ an everywhere non-zero vector field
$\omega$ defined by the formula $\omega_{\zeta}=\zeta$, i\. e\. assigning to
a point $\zeta\in U_0$ the vector $\zeta$. This vector field $\omega$
lies in the distribution $\Cal D_2$. It is easy to see that the
1-dimensional distribution $\Cal I$ generated by the vector field
$\omega$ and the 1-dimensional distribution $\Cal D_1$ are transversal.

Fixing a volume form $\theta_0\in\Lambda^6V^*$, we get for each $\zeta\in
U_0$ a tangent structure $F(\zeta)$. Namely, this tangent structure can
be determined by the formula
$$
(\iota_v\zeta)\wedge\zeta=\iota_{F(\zeta)v}\theta_0.
$$
For any 3-form $\Omega\in\Lambda^3V^*$ we can then define
$$
(D_{F(\zeta)}\Omega)(v_1,v_2,v_3)=\Omega(F(\zeta)v_1,v_2,v_3)+
\Omega(v_1,F(\zeta)v_2,v_3)+\Omega(v_1,v_2,F(\zeta)v_3).
$$
It is obvious that if $\Omega\in T_{\zeta}U_0$, then also $D_F\Omega
\in T_{\zeta}U_0$. Consequently, on $T_{\zeta}U_0$ we can define
an endomorphism $\Cal N(\zeta)$ by the formula $\Cal N(\zeta)=D_{F(\zeta)}$.
In this way we get on $U_0$ a tensor field $\Cal N$ of type $(1,1)$. It is
easy to see that $\Cal N^3=0$.

Our main aim in this section will be to prove the following proposition.
\proclaim{6.3. Proposition}
On $U_0$ we have the following chain of distributions:
$$
\im\Cal N^2\subset\ker\Cal N\subset\im\Cal N\subset\ker\Cal N^2,
$$
where $\im\Cal N^2=\Cal D_1$ and $\im\Cal N=\Cal D_2$. The distributions
$\im Cal N^2$, $\ker\Cal N$, and $\im\Cal N$ are integrable. The
distribution $\ker\Cal N^2$ is not integrable.
\endproclaim

\definition{6.4. Remark}
If $A\in End(V)$ is arbitrary we can define $D_A\Omega$ for any $\Omega\in
\Lambda^kV^*$ by the formula
$$
(D_F\Omega)(v_1,\dots,v_k)=\sum_{i=1}^k\Omega(v_1,\dots,v_{i-1},Av_i,
v_{i+1},\dots,v_k).
$$
It is well known that $D_A$ is a derivation on the graded algebra
$\Lambda^*V^*$.
\enddefinition

We shall first investigate the subspace $\im\Cal N^2$. Let $\Omega\in
\im\Cal N^2(\zeta)$. If $\Omega=\Cal N^2(\zeta)\tilde{\Omega}$, then we
have
$$
\Omega(v_1,v_2,v_3)=2(\tilde{\Omega}(Fv_1,Fv_2,v_3)+
\tilde{\Omega}(Fv_1,v_2,Fv_3)+\tilde{\Omega}(v_1,Fv_2,Fv_3)),
$$
where $F=F(\zeta)$. It is easy to see that if one of the entries
$v_1,v_2,v_3$ belongs to $V_0(\zeta)$, then $\Omega(v_1,v_2,v_3)=0$, or in
other words, $\Omega\in\Cal D_1(\zeta)$. Because obviously $\im\Cal N^2
\ne0$, we get easily the following lemma. (Notice that $\dim\im\Cal
N^2=1$.)
\proclaim{6.5. Proposition}
There is $\im\Cal N^2=\Cal D_1$ and $\im\Cal N^2\subset\ker\Cal N$. The
distribution $\im\Cal N^2\subset TU_0$ is integrable.
\endproclaim

Next, we shall consider the subspace $\Cal D_2(\zeta)$. It is obvious
that for any $\Omega\in\Cal D_2(\zeta)$ the correspondence $v\in
V_0(\zeta)\mapsto\iota_v\Omega$ defines a homomorphism
$$
\kappa_{\Omega}:V_0(\zeta)\rightarrow\Lambda^2W(\zeta)^*.
$$
We have obvious formulas
$$
\kappa_{\Omega+\tilde{\Omega}}=\kappa_{\Omega}+\kappa_{\tilde{\Omega}},
\quad\kappa_{a\Omega}=a\kappa_{\Omega}
$$
for any $\Omega,\tilde{\Omega}\in\Cal D_2(\zeta)$ and any $a\in\Bbb R$.
Using the isomorphism $\kappa_{\zeta}:V_0(\zeta)\rightarrow\Lambda^2
W(\zeta)^*$, we can define a homomorphism
$$
k_{\Omega}:\Cal D_2(\zeta)\rightarrow End(V_0(\zeta)),\quad k_{\Omega}(v)=
\kappa_{\zeta}^{-1}\kappa_{\Omega}(v).
$$
It is easy to see that $\ker k_{\Omega}=\Cal D_1(\zeta)$. Consequently,
we get a monomorphism
$$
K_{\Omega}:\Cal D_2(\zeta)/\Cal D_1(\zeta)\rightarrow End(V_0(\zeta)).
$$
Because $\dim\Cal D_2(\zeta)/\Cal D_1(\zeta)=\dim End(V_0(\zeta))=9$, we
can see that $K_{\Omega}$ is an isomorphism.

\proclaim{6.6. Proposition}
There is $\im\Cal N=\Cal D_2$ and $\dim\im\Cal N=10$.
\endproclaim
\demo{Proof}
If $\Omega=\Cal N(\zeta)\hat{\Omega}$, where $\hat{\Omega}\in
T_{\zeta}U_0$, we have
$$
\Omega(v_1,v_2,v_3)=\hat{\Omega}(F(\zeta)v_1,v_2,v_3)+
                    \hat{\Omega}(v_1,F(\zeta)v_2,v_3)+
                    \hat{\Omega}(v_1,v_2,F(\zeta)v_3),
$$
and it is obvious that $\Omega\in\Cal D_2(\zeta)$. This shows that $\im
\Cal N\subset\Cal D_2$.

Conversely, let us assume that $\Omega\in\Cal D_2(\zeta)$. We choose a
basis $v_1,v_2,v_3$ of $V_0(\zeta)$, and we denote $\pi(\zeta):
V\rightarrow W(\zeta)$ the projection. Because $\Omega\in\Cal
D_2(\zeta)$, there exist 2-forms $\tilde{\Omega}_1,\tilde{\Omega}_2,
\tilde{\Omega}_3\in\Lambda^2W(\zeta)^*$ such that
$$
\iota_{v_i}\Omega=\pi(\zeta)^*\tilde{\Omega}_i,\quad i=1,2,3.
$$
Let us take now 1-forms $\beta_1,\beta_2,\beta_3\in V^*$ such that
$\beta_i(v_j)=\delta_{ij}$. We shall consider a 3-form
$$
\hat{\Omega}=\sum_{i=1}^3\beta_i\wedge\pi(\zeta)^*\tilde{\Omega}_i.
$$
Now we can easily see that $\iota_v(\Omega-\hat{\Omega})=0$ for any
$v\in V_0(\zeta)$, or in other words $\Omega-\hat{\Omega}\in\Cal D_1=
\im\Cal N^2$. This means that there is a 3-form $\bar{\Omega}\in T_{\zeta}
U_0$ such that $\Omega-\hat{\Omega}=\Cal N^2(\zeta)\bar{\Omega}$.

Let us consider the monomorphism $\pi(\zeta)^*:\Lambda^*W(\zeta)^*
\rightarrow\Lambda^*V^*$. It is easy to see that $\pi(\zeta)^*
W(\zeta)^*$ has a basis $D_{F(\zeta)}\beta_1,D_{F(\zeta)}\beta_2,
D_{F(\zeta)}\beta_3$, and that
$$
D_{F(\zeta)}^2\beta_1=D_{F(\zeta)}^2\beta_2=D_{F(\zeta)}^2\beta_3=0.
$$
It is obvious that any 2-form $\Omega'\in\pi(\zeta)^*\Lambda^2
W(\zeta)^*$ belongs to $\im D_{F(\zeta)}^2$. Consequently, we can find
2-forms $\Omega'_1,\Omega'_2,\Omega'_3$ such that
$$
\pi(\zeta)^*\tilde{\Omega}_i=D_{F(\zeta)}^2\Omega'_i.
$$
We have then
$$
\gather
\hat{\Omega}=\sum_{i=1}^3\beta_i\wedge\pi(\zeta)^*\tilde{\Omega}_i=
\sum_{i=1}^3\beta_i\wedge D_{F(\zeta)}^2\Omega'_i=\\
=\sum_{i=1}^3D_{F(\zeta)}(\beta_i\wedge D_{F(\zeta)}\Omega'_i)
-\sum_{i=1}^3D_{F(\zeta)}\beta_i\wedge D_{F(\zeta)}\Omega'_i=\\
=D_{F(\zeta)}\sum_{i=1}^3\beta_i\wedge D_{F(\zeta)}\Omega'_i
-\sum_{i=1}^3D_{F(\zeta)}(D_{F(\zeta)}\beta_i\wedge\Omega'_i)=\\
=D_{F(\zeta)}\sum_{i=1}^3(\beta_i\wedge D_{F(\zeta)}\Omega'_i
-D_{F(\zeta)}\beta_i\wedge\Omega'_i).
\endgather
$$
Now we can see that $\Omega\in\im\Cal N(\zeta)$, which finishes the
proof.
\enddemo

\proclaim{6.7. Proposition}
There is the inclusion $\ker\Cal N\subset\im\Cal N$.
\endproclaim
\demo{Proof}
Let $\Omega\in\ker\Cal N(\zeta)$. Then we have (we write $F$ instead of
$F(\zeta)$)
$$
\Omega(Fv_1,v_2,v_3)+\Omega(v_1,Fv_2,v_3)+\Omega(v_1,v_2,Fv_3)=0.
$$
Using this relation we get
$$
\align
&\Omega(Fv_1,Fv_2,v_3)=-\Omega(v_1,F^2v_2,v_3)-\Omega(v_1,Fv_2,Fv_3)=
-\Omega(v_1,Fv_2,Fv_3)\\
&\Omega(Fv_1,Fv_2,v_3)=-\Omega(F^2v_1,v_2,v_3)-\Omega(Fv_1,v_2,Fv_3)=
-\Omega(Fv_1,v_2,Fv_3)
\endalign
$$
Adding these two relations, we obtain
$$
\gather
2\Omega(Fv_1,Fv_2,v_3)=-\Omega(v_1,Fv_2,Fv_3)-\Omega(Fv_1,v_2,Fv_3),\\
\Omega(Fv_1,Fv_2,v_3)=-\Omega(Fv_1,Fv_2,v_3)-\Omega(Fv_1,v_2,Fv_3)
-\Omega(v_1,Fv_2,Fv_3)=\\
=-\frac{1}{2}D_F^2\Omega(v_1,v_2,v_3)=0,
\endgather
$$
which shows that $\Omega\in\Cal D_2(\zeta)$.
\enddemo

\proclaim{6.8. Proposition}
Let $\zeta\in U_0$. Then $\Omega\in T_{\zeta}U_0$ belongs to $\ker\Cal
N^2$ if and only if $\zeta\wedge\Omega=0$. Moreover $\dim\ker\Cal N^2=
18$.
\endproclaim
\demo{Proof}
Let us choose vectors $v,v',v''\in V$ such that $Fv,Fv',Fv'',v,v',v''$ is
a basis of $V$. (We denote for simplicity $F=F(\zeta)$.) We shall
consider the value $(\zeta\wedge\Omega)(Fv,Fv',Fv'',v,v',v'')$. (We
recall that $\zeta(w,w',\cdot)=0$ if $w,w'\in V_0(\zeta)$,
$\zeta(w,Fw,\cdot)=0$ for any $w\in V$, and $\Omega|V_0(\zeta)=0$.) We get
$$
\gather
(\zeta\wedge\Omega)(Fv,Fv',Fv'',v,v',v'')=
\zeta(Fv,v',v'')\Omega(Fv',Fv'',v)+\\
+\zeta(Fv',v,v'')\Omega(Fv,Fv'',v')+
\zeta(Fv'',v,v')\Omega(Fv,Fv',v'')=\\
=\zeta(Fv,v',v'')[\Omega(Fv,Fv',v'')+\Omega(Fv,v',Fv'')+
                   \Omega(v,Fv',Fv'')].
\endgather
$$
Because $\zeta(Fv,v',v'')\ne0$ the first assertion easily follows. Now
it is obvious that $\dim\ker\Cal N^2=18$.
\enddemo

\proclaim{6.9. Proposition}
There is $\im\Cal N\subset\ker\Cal N^2$.
\endproclaim
\demo{Proof}
If $\Omega\in\im\Cal N(\zeta)=\Cal D_2(\zeta)$ then obviously $\zeta
\wedge\Omega=0$.
\enddemo

On the trivial vector bundle $\Cal V$ with fiber $V$ over $U_0$ we
introduce a linear connection $\nabla$. For any vector field $\Omega$ on
$U_0$ and any section $S$ of $\Cal V$ we define $\nabla_{\Omega}S=\Omega S$.
Obviously, $\nabla$ induces a linear connection on every exterior power
$\Lambda^kV^*$, which will be denoted by the same symbol. It is obvious
that the same formula $\bar{\nabla}_{\bar{\Omega}}\bar{S}=\bar{\Omega}\bar{S}$,
where $\bar{S}$ is a section of the trivial vector bundle $\bar{\Cal V}$
with fiber $V$ over $\Lambda^3V^*$, extends the connection $\nabla$ to
the whole vector space $\Lambda^3V^*$. The connection $\bar{\nabla}$ induces
again a linear connection on the vector bundle $\Lambda^k\bar{\Cal V}^*$,
which will be denoted again by the symbol $\bar{\nabla}$. Let $\Omega_1$
and $\Omega_2$ be (local) vector fields on $U_0$, and let
$\bar{\Omega}_1$ and $\bar{\Omega}_2$ be their (local)
extensions. Because the connection $\bar{\nabla}$ is flat, we have
$\bar{\nabla}_{\bar{\Omega}_1}\bar{\Omega_2}-\bar{\nabla}_{\bar{\Omega}_2}
\bar{\Omega}_1=[\bar{\Omega}_1,\bar{\Omega}_2]$. Restricting this formula
to the submanifold $U_0$, we obtain the formula
$$
\nabla_{\Omega_1}\Omega_2-\nabla_{\Omega_2}\Omega_1=[\Omega_1,\Omega_2],
$$
which will be needed in the sequel.

\proclaim{6.10. Lemma}
Let $S$ be a section of the subbundle $\Cal V_0$, and let $\Omega$ be a
vector field on $U_0$ lying in $\im\Cal N$. Then $\nabla_{\Omega}S$ is
also a section of the subbundle $\Cal V_0$.
\endproclaim
\demo{Proof}
Because $S$ is a section of the subbundle $\Cal V_0$, we have the
relation $(\iota_S\omega)\wedge\omega=0$. Applying to this relation
$\nabla_{\Omega}$, we obtain
$$
(\iota_{\nabla_{\Omega}S}\omega)\wedge\omega+(\iota_S\Omega)\wedge\omega+
(\iota_S\omega)\wedge\Omega=0.
$$
It is easy to see that the second term vanishes. The last term vanishes
by virtue of Lemma 6.3. Consequently, we obtain
$(\iota_{\nabla_{\Omega}S}\omega)\wedge\omega=0$, which shows that
$\nabla_{\Omega}S$ is a section of $\Cal V_0$.
\enddemo

\definition{6.11. Remark}
The previous lemma shows that the connection $\nabla$ on $\Cal V$
induces a partial connection on $\Cal V_0$, which we shall denote by the
same symbol. This partial connection determines the covariant derivative
$\nabla_{\Omega}$ only for vector the fields $\Omega$ lying in $\im\Cal
N$. This partial connection induces a partial connection on the vector
bundle $\Cal W$ and on any exterior power of the vector bundles $\Cal
V_0$ and $\Cal W$. Moreover, if $\tilde{\Omega}$ is a vector field on
$U_0$ (i\. e\. a section of $\Lambda^3\Cal V^*$ such that
$\tilde{\Omega}|\Cal V_0=0$), then for any vector field $\Omega$ lying in
$\im\Cal N$ and any three sections $S_1,S_2,S_3$ of $\Cal V_0$ we have
$$
\gather
\tilde{\Omega}(S_1,S_2,S_3)=0\\
\nabla_{\Omega}(\tilde{\Omega}(S_1,S_2,S_3))=0\\
(\nabla_{\omega}\tilde{\Omega})(S_1,S_2,S_3)+
\tilde{\Omega}(\nabla_{\Omega}S_1,S_2,S_3)+
\tilde{\Omega}(S_1,\nabla_{\Omega}S_2,S_3)+
\tilde{\Omega}(S_1,S_2,\nabla_{\Omega}S_3)=0\\
(\nabla_{\Omega}\tilde{\Omega})(S_1,S_2,S_3)=0,
\endgather
$$
which shows that the partial connection $\nabla$ induces a partial
connection (again denoted by the same symbol) on $TU_0$. Because the
original connection on $\Cal V$ is flat, we have for any two vector
fields $\Omega$ and $\tilde{\Omega}$ lying in $\im\Cal N$
$$
\nabla_{\Omega}\tilde{\Omega}-\nabla_{\tilde{\Omega}}\Omega=
[\Omega,\tilde{\Omega}].
$$
\enddefinition

\proclaim{6.12. Proposition}
The distribution $\im\Cal N$ is integrable.
\endproclaim
\demo{Proof}
According to Proposition 6.6 there is $\im\Cal N=\Cal D_2$. Let us take
two vector fields $\Omega,\tilde{\Omega}$ lying in $\Cal D_2$, and three
sections $S_1,S_2,S_3$ of $\Cal V$ such that $S_1$ and $S_2$ lie in
$\Cal V_0$. Then we have
$$
\gather
(\nabla_{\Omega}\tilde{\Omega})(S_1,S_2,S_3)=
\nabla_{\Omega}(\tilde{\Omega}(S_1,S_2,S_3))-\\
-\tilde{\Omega}(\nabla_{\Omega}S_1,S_2,S_3)
-\tilde{\Omega}(S_1,\nabla_{\Omega}S_2,S_3)
-\tilde{\Omega}(S_1,S_2,\nabla_{\Omega}S_3)=0\\
\endgather
$$
according to Lemma 6.10. This shows that $\nabla_{\Omega}\tilde{\Omega}$
lies in $\Cal D_2$. Now, it is obvious that $[\Omega,\tilde{\Omega}]=
\nabla_{\Omega}\tilde{\Omega}-\nabla_{\tilde{\Omega}}\Omega$ lies in
$\Cal D_2$.
\enddemo

\proclaim{6.13. Proposition}
$\ker\Cal N=\{\Omega\in\im\Cal N;\Tr k(\Omega)=0\}$ and $\dim
\ker\Cal N=9$.
\endproclaim
\demo{Proof}
We shall denote for simplicity $V_0=V_0(\zeta)$, $F=F(\zeta)$,
$W=W(\zeta)$, $\pi=\pi(\zeta)$. Let us notice first that for each
endomorphism $A\in End(V_0)$ there exists an endomorphism $B\in End(V)$
(not uniquely determined) such that
$$
AF=FB\quad\text{and}\quad BV_0\subset V_0.
$$
Moreover, any endomorphism $B$ with these properties induces
an endomorphism $\tilde{B}\in End(W)$ and $\Tr\tilde{B}=\Tr A$.

Let us take now a 3-form $\Omega\in\im\Cal N(\zeta)=\Cal D_2(\zeta)$. We
have
$$
(\Cal N(\zeta)\Omega)(v,v',v'')=\Omega(Fv,v',v'')+\Omega(v,Fv',v'')+
\Omega(v,v',Fv'').
$$
It is easy to see that $\Cal N(\zeta)\Omega\in\Cal D_1(\zeta)$, and
consequently there exists a uniquely determined 3-form $\tilde{\Omega}
\in\Lambda^3W^*$ such that $\Cal N(\zeta)\Omega=\pi^*\tilde{\Omega}$.
Similarly, there is a 3-form $\tilde{\zeta}\in\Lambda^3W^*$ such that
$\Cal N(\zeta)\zeta=\pi^*\tilde{\zeta}$. We recall that the homomorphism
$\pi^*:\Lambda^3W^*\rightarrow\Lambda^3V^*$ is a monomorphism.
Consequently $\tilde{\zeta}\ne0$.

Let us take now $A=k(\zeta)$. Obviously for any $v,v',v''\in V$ we have
$$
\align
\zeta(AFv,v',v'')&=\Omega(Fv,v',v''),\\
\zeta(v,AFv',v'')&=\Omega(v,Fv',v''),\\
\zeta(v,v',AFv'')&=\Omega(v,v',Fv'').
\endalign
$$
Then we get
$$
\gather
(\Cal N(\zeta)\Omega)(v,v',v'')=\Omega(Fv,v',v'')+\Omega(v,Fv',v'')+
\Omega(v,v',Fv'')=\\
=\zeta(AFv,v',v'')+\zeta(v,AFv',v'')+\zeta(v,v',AFv'')=\\
=\zeta(FBv,v',v'')+\zeta(v,FBv',v'')\zeta(v,v',FBv'')=\\
=(1/3)[\zeta(FBv,v',v'')+\zeta(Bv,Fv',v'')+\zeta(Bv,v',Fv'')]+\\
+(1/3)[\zeta(Fv,Bv',v'')+\zeta(v,FBv',v'')+\zeta(v,Bv',Fv'')]+\\
+(1/3)[\zeta(Fv,v',Bv'')+\zeta(v,Fv',Bv'')+\zeta(v,v',FBv'')]=\\
=(1/3)\tilde{\zeta}(\tilde{B}[v],[v'],[v''])
+(1/3)\tilde{\zeta}([v],\tilde{B}[v'],[v''])
+(1/3)\tilde{\zeta}([v],[v'],\tilde{B}[v''])=\\
=(1/3)\Tr(\tilde{B})\tilde{\zeta}([v],[v'],[v''])
=(1/3)\Tr(A)\zeta(v,v',v''),
\endgather
$$
which shows that $\Cal N(\zeta)\Omega=0$ if and only if $\Tr(A)=0$.
\enddemo

\proclaim{6.14. Lemma}
Let $M$ be a differentiable manifold, and let $\xi$ be an
$n$-dimensional differentiable vector bundle over $M$ endowed with a
linear connection $\nabla$. Let $A$ be an endomorphism of the vector
bundle $\xi$, i\. e\. a section of the vector bundle $\xi^*\otimes\xi$.
Then for any vector field $X$ on $M$ we have
$$
\Tr(\nabla_XA)=X\Tr(A).
$$
\endproclaim
\demo{Proof}
Let us choose (at least locally) a non-zero $n$-form $\varepsilon$ on
$\xi$. Then for any vector fields $X_1,\dots,X_n$ we have
$$
\sum_{i=1}^n\varepsilon(X_1,\dots,X_{i-1},AX_i,X_{i+1},\dots,X_n)=\Tr(A)
\cdot\varepsilon(X_1,\dots,X_n).
$$
Let $X$ be a vector field on $M$. Applying $\nabla_X$ to the above
equality, we obtain
$$
\sum_{i=1}^n\varepsilon(X_1,\dots,X_{i-1},(\nabla_XA)X_i,X_{i+1},\dots,
X_n)=(X\Tr(A))\cdot\varepsilon(X_1,\dots,X_n),
$$
which implies the desired equality.
\enddemo

\proclaim{6.15. Proposition}
The distribution $\ker\Cal N$ is integrable.
\endproclaim
\demo{Proof}
Let $\Omega$ and $\tilde{\Omega}$ be two vector fields lying in the
distribution $\ker\Cal N$. We denote $A=k_{\Omega}$ and $\tilde{A}=
k_{\tilde{\Omega}}$. According to the previous result there is $\Tr(A)=
\Tr(\tilde{A})=0$. For any section $S$ of $\Cal V_0$ and any constant
sections $S',S''$ of $\Cal V$ we have
$$
\omega(AS,S',S'')=\Omega(S,S',S''),\quad\omega(\tilde{A}S,S',S'')=
\tilde{\Omega}(S,S',S'').
$$
Applying $\nabla_{\Omega}$ to the second equality we obtain
$$
\gather
(\nabla_{\Omega}\omega)(\tilde{A}S,S',S'')+
\omega((\nabla_{\Omega}\tilde{A})S,S',S'')+
\omega(\tilde{A}\nabla_{\Omega}S,S',S'')=\\
=(\nabla_{\Omega}\tilde{\Omega})(S,S',S'')+
\tilde{\Omega}(\nabla_{\Omega}S,S',S'')\\
\Omega(\tilde{A}S,S',S'')+
\omega((\nabla_{\Omega}\tilde{A})S,S',S'')=
(\nabla_{\Omega}\tilde{\Omega})(S,S',S'')\\
\omega(A\tilde{A}S,S'.S'')+
\omega((\nabla_{\Omega}\tilde{A})S,S',S'')=
(\nabla_{\Omega}\tilde{\Omega})(S,S',S'')\\
(\nabla_{\Omega}\tilde{\Omega})(S,S',S'')=
\omega((A\tilde{A}+\nabla_{\Omega}\tilde{A})S,S',S'').
\endgather
$$
Similarly we obtain
$$
(\nabla_{\tilde{\Omega}}\Omega)(S,S',S'')=
\omega((\tilde{A}A+\nabla_{\tilde{\Omega}}A)S,S',S'').
$$
Substracting the last two equalities we have
$$
[\Omega,\tilde{\Omega}](S,S',S'')=
(\nabla_{\Omega}\tilde{\Omega}-\nabla_{\tilde{\Omega}}\Omega)
(S,S',S'')=
\omega(([A,\tilde{A}]+
\nabla_{\Omega}\tilde{A}-\nabla_{\tilde{\Omega}}A)S,S',S''),
$$
which shows that
$$
k_{[\Omega,\tilde{\Omega}]}=[A,\tilde{A}]+
\nabla_{\Omega}\tilde{A}-\nabla_{\tilde{\Omega}}A.
$$
On any integral submanifold of the distribution $\im\Cal N$ we have
$$
\Tr([A,\tilde{A}]+
\nabla_{\Omega}\tilde{A}-\nabla_{\tilde{\Omega}}A)=
0+\Omega\Tr(\tilde{A})-\tilde{\Omega}\Tr(A)=0.
$$
This finishes the proof.
\enddemo

\proclaim{6.16. Proposition}
The distribution $\ker\Cal N^2$ is not integrable.
\endproclaim
\demo{Proof}
Let $\Omega$ and $\tilde{\Omega}$ be two vector fields on $U_0$ lying in $\ker
\Cal N^2$. We shall apply the vector field $\Omega$ to the relation
$\tilde{\Omega}\wedge\omega=0$. We get
$$
(\nabla_{\Omega}\tilde{\Omega})\wedge\omega+\tilde{\Omega}\wedge\Omega=0.
$$
Interchanging $\Omega$ and $\tilde{\Omega}$ and substracting the two
relations, we obtain
$$
\gather
[\Omega,\tilde{\Omega}]\wedge\omega+
\tilde{\Omega}\wedge\Omega-\Omega\wedge\tilde{\Omega}=0\\
[\Omega,\tilde{\Omega}]\wedge\omega=2\Omega\wedge\tilde{\Omega}.
\endgather
$$
Let us choose now vectors $\alpha,\tilde{\alpha}\in T_{\omega_0}U_0$ as
follows:
$$
\alpha=\alpha_1\wedge\alpha_2\wedge\alpha_5,\quad
\tilde{\alpha}=\alpha_3\wedge\alpha_4\wedge\alpha_6.
$$
It is easy to verify that $\Omega,\tilde{\Omega}\in\ker\Cal
N^2(\omega_0)$. We choose now vector fields $\Omega$ and $\tilde{\Omega}$
in such a way that they lie in $\ker\Cal N^2$ and
$\Omega_{\omega_0}=\alpha$ and $\tilde{\Omega}_{\omega_0}=\tilde{\alpha}$.
According to the above formula we have then
$$
[\Omega,\tilde{\Omega}]_{\omega_0}\wedge\omega_0=
2\alpha\wedge\tilde{\alpha}\ne0,
$$
which shows that the vector field $[\Omega,\tilde{\Omega}]$ does not lie
in $\ker\Cal N^2$.
\enddemo

\enddocument
\end
*********************************************

Let $V$ be a 6-dimensional real vector space, and let $\omega$ be a
multisymplectic 3-form on $V$ such that $V_0=\Delta^2(\omega)$ is a
3-dimensional subspace. Let us consider the injective homomorphism
$$
\kappa:V\rightarrow\Lambda^2V^*,\quad\kappa v=\iota_v\omega.
$$
The image $\kappa V_0$, which is a 3-dimensional subspace in
$\Lambda^2V^*$, consists of decomposable 2-forms. Then according to [V]
there are two possibilities. Either, there exist linerly independent
1-forms $\beta_1,\beta_2,\beta_3\in V^*$ such that
$$
\beta_1\wedge\beta_2,\quad\beta_1\wedge\beta_3,\quad\beta_2\wedge\beta_3
$$
is a basis of $\kappa V_0$. Or, there exist linearly independent 1-forms
$\beta_0,\beta_1,\beta_2,\beta_3\in V^*$ such that
$$
\beta_0\wedge\beta_1,\quad\beta_0\wedge\beta_2,\quad\beta_0\wedge\beta_3
$$
is a basis of $\kappa V_0$.

We shall start to consider the second possibility. Let $f_1,f_2,f_3$ be
the basis of $V_0$ such that
$$
\kappa f_1=\beta_0\wedge\beta_1,\quad\kappa f_2=\beta_0\wedge\beta_2,
\quad\kappa f_3=\beta_0\wedge\beta_3.
$$
We find easily that
$$
0=\iota_{f_1}\iota_{f_1}\omega=\iota_{f_1}\kappa f_1=\iota_{f_1}
(\beta_0\wedge\beta_1)=\beta_0(f_1)\beta_1-\beta_1(f_1)\beta_0,
$$
which shows that $\beta_0(f_1)=0$. Along the same lines we find that
$\beta_0(f_2)=0$ and $\beta_0(f_3)=0$. In other words, we have shown
that $\beta_0|V_0=0$. We choose now a basis $\alpha_1,\dots,\alpha_6$ of
$V^*$ such that $\alpha_4|V_0=\alpha_5|V_0=\alpha_6|V_0=0$ and $\alpha_6
=\beta_0$. We can express $\omega$ in the form
$$
\omega=\eta\wedge\alpha_6+\zeta,
$$
where $\eta$ and $\zeta$ do not contain $\alpha_6$. We have
$$
\iota_{f_1}\omega=(\iota_{f_1}\eta)\wedge\alpha_6+\iota_{f_1}\zeta,
$$
which shows that $\iota_{f_1}\zeta=0$. Similarly we find that
$\iota_{f_2}\zeta=0$ and $\iota_{f_3}\zeta=0$. Summarizing, we can see
that the 3-form $\zeta$ does not contain neither $\alpha_6$, nor
$\alpha_1$, $\alpha_2$, and $\alpha_3$. Therefore we have $\zeta=0$.
We denote $N=\{v\in V;\alpha_6(v)=0\}$. $N$ is a 5-dimensional (hence odd
dimensional) subspace, and consequently there exists a nonzero vector
$w\in N$ such that $\iota_w(\eta|N)=0$. We get then
$$
\iota_w\omega=\iota_w(\eta\wedge\alpha_6)=(\iota_w\eta)\wedge\alpha_6+
\alpha_6(w)\eta=0.
$$
This shows that the 3-form $\omega$ is not multisymplectic, which is a
contradiction.

This means that there exist 1-forms $\beta_1,\beta_2,\beta_3$ such that
$\beta_2\wedge\beta_3,\beta_3\wedge\beta_1,\beta_1\wedge\beta_2
$ is a basis of $\kappa V_0$. We take $f_1,f_2,f_3\in V_0$ such that
$$
\kappa f_1=\beta_2\wedge\beta_3,\quad
\kappa f_2=\beta_3\wedge\beta_1,\quad
\kappa f_3=\beta_1\wedge\beta_2.
$$
This time we get
$$
0=\iota_{f_1}\iota_{f_1}\omega=\iota_{f_1}(\beta_2\wedge\beta_3)=
\beta_2(f_1)\beta_3-\beta_3(f_1)\beta_2,
$$
which implies $\beta_2(f_1)=0$ and $beta_3(f_1)=0$. Proceeding in this
way, we obtain finally
$$
\beta_2(f_1)=\beta_3(f_1)=0,\quad
\beta_1(f_2)=\beta_3(f_2)=0,\quad
\beta_1(f_3)=\beta_2(f_3)=0.
$$
Further we get
$$
\align
\iota_{f_1}\iota_{f_2}\omega&=\iota_{f_1}(\beta_3\wedge\beta_1)=
\beta_3(f_1)\beta_1-\beta_1(f_1)\beta_3=-\beta_1(f_1)\beta_3,\\
-\iota_{f_2}\iota_{f_1}\omega&=-\iota_{f_2}(\beta_2\wedge\beta_3)=
\beta_3(f_2)\beta_2-\beta_2(f_2)\beta_3=-\beta_2(f_2)\beta_3,
\endalign
$$
which implies $\beta_1(f_1)=\beta_2(f_2)$. Proceeding further, we find
easily that
$$
\beta_1(f_1)=\beta_2(f_2)=\beta_3(f_3).
$$
Now we must consider two possibilities. The first possibility is that
$\beta_1(f_1)=\beta_2(f_2)=\beta_3(f_3)=a\ne0$. Obviously, we can
complete the 1-forms $\beta_1,\beta_2,\beta_3$ to a basis $\beta_1,
\beta_2,\beta_3,\beta_4,\beta_5,\beta_6$ in such a way that
$\beta_4|V_0=\beta_5|V_0=\beta_6|V_0=0$. Now we can express the 3-form
$\omega$ in the following way.
$$
\omega=\eta_1\wedge\beta_1+\eta_2\wedge\beta_2+\eta_3\wedge\beta_3+\zeta,
$$
where $\eta_1$ does not contain $\beta_1$, $\eta_2$ contains neither
$\beta_1$ nor $\beta_2$, and $\eta_3$ and $\zeta$ contain neither
$\beta_1$ nor $\beta_2$ nor $\beta_3$. We can see that
$$
\zeta=b\beta_4\wedge\beta_5\wedge\beta_5.
$$
We find easily
$$
\beta_2\wedge\beta_3=\iota_{f_1}\omega=(\iota_{f_1}\eta_1)\wedge\beta_1+
a\eta_1,
$$
$$
\omega=(1/a)\beta_1\wedge\beta_2\wedge\beta_3+\eta_2\wedge\beta_2+
\eta_3\wedge\beta_3+b\beta_4\wedge\beta_5\wedge\beta_5,
$$
$$
\beta_3\wedge\beta_1=\iota_{f_2}\omega=-\beta_1\wedge\beta_3+a\eta_2.
$$
This last equation shows that $\eta_2=0$. Thus we get
$$
\omega=(1/a)\beta_1\wedge\beta_2\wedge\beta_3+
\eta_3\wedge\beta_3+b\beta_4\wedge\beta_5\wedge\beta_5,
$$
$$
\beta_1\wedge\beta_2=\iota_{f_3}\omega=\beta_1\wedge\beta_2+
a\eta_3,
$$
which shows this time that $\eta_3=0$. Finally we get
$$
\omega=(1/a)\beta_1\wedge\beta_2\wedge\beta_3+
b\beta_4\wedge\beta_5\wedge\beta_5.
$$
Here we come again to a contradiction. If $b=0$ the 3-form $\omega$ is
not multisymplectic. If $b\ne0$ then $\Delta^2(\omega)$ is the union of
$V_0$ with a transversal 3-dimensional subspace. This shows that only
the second possibility can occur, namely that we have
$$
\beta_1(f_1)=\beta_2(f_2)=\beta_3(f_3)=0.
$$

This means that we shall consider the second possibility. Let us change
first our notation. We set $\alpha_4=\beta_1$, $\alpha_5=
\beta_2$, and $\alpha_6=\beta_3$. We complete the 1-forms $\alpha_4,
\alpha_5,\alpha_6$ to a basis $\alpha_1,\alpha_2,\alpha_3,\alpha_4,
\alpha_5,\alpha_6$ in such a way that $f_1,f_2,f_3$ and $\alpha_1|V_0,
\alpha_2|V_0,\alpha_3|V_0$ are dual bases. This time we shall write
$$
\omega=\eta_1\wedge\alpha_1+\eta_2\wedge\alpha_2+\eta_3\wedge\beta_3+
\zeta,
$$
where $\eta_1$ does not contain $\alpha_1$, $\eta_2$ contain neither
$\alpha_1$ nor $\alpha_2$, and $\eta_3$ and $\zeta$ contain neither
$\alpha_1$ nor $\alpha_2$ nor $\alpha_3$. We have obviously
$$
\zeta=b\alpha_4\wedge\alpha_5\wedge\alpha_6,
$$
$$
\alpha_5\wedge\alpha_6=\iota_{f_1}\omega=\eta_1,
$$
$$
\omega=\alpha_1\wedge\alpha_5\wedge\alpha_6+\eta_2\wedge\alpha_2+
\eta_3\wedge\beta_3+\zeta,
$$
$$
\alpha_6\wedge\alpha_4=\iota_{f_2}\omega=\eta_2,
$$
$$
\omega=\alpha_1\wedge\alpha_5\wedge\alpha_6+
       \alpha_2\wedge\alpha_6\wedge\alpha_4+
       \eta_3\wedge\alpha_3+\zeta
$$
$$
\alpha_4\wedge\alpha_5=\iota_{f_3}\omega=\eta_3,
$$
$$
\omega=\alpha_1\wedge\alpha_5\wedge\alpha_6+
       \alpha_2\wedge\alpha_6\wedge\alpha_4+
       \alpha_3\wedge\alpha_4\wedge\alpha_5+
       b\alpha_4\wedge\alpha_5\wedge\alpha_6
$$

We choose now arbitrarily a nonzero $\theta\in\Lambda^6V^*$. Then for
any vector $v\in V$ there exists a unique vector $Qv\in V$ such that
$$
(\iota_v\omega)\wedge\omega=\iota_{Qv}\theta.
$$

Let $f_1,\dots,f_6$ be the dual basis to $\alpha-1,\dots,\alpha_6$. We
take $v=c_1f_1+\dots+c_6f_6$, Then we get
$$
\align
\iota_v\omega=&c_1\alpha_5\wedge\alpha_6-
               c_5\alpha_1\wedge\alpha_6+
               c_6\alpha_1\wedge\alpha_5\\
             +&c_2\alpha_6\wedge\alpha_4-
               c_6\alpha_2\wedge\alpha_4+
               c_4\alpha_2\wedge\alpha_6\\
             +&c_3\alpha_4\wedge\alpha_5-
               c_4\alpha_3\wedge\alpha_5+
               c_5\alpha_3\wedge\alpha_4\\
             +&bc_4\alpha_5\wedge\alpha_6-
               bc_5\alpha_4\wedge\alpha_6+
               bc_6\alpha_4\wedge\alpha_5\\
\endalign
$$
$$
(\iota_v\omega)\wedge\omega=
-2c_4\alpha_2\wedge\alpha_3\wedge\alpha_4\wedge\alpha_5\wedge\alpha_6
+2c_5\alpha_1\wedge\alpha_3\wedge\alpha_4\wedge\alpha_5\wedge\alpha_6
-2c_6\alpha_1\wedge\alpha_2\wedge\alpha_4\wedge\alpha_5\wedge\alpha_6.
$$
Now we get easily the following lemma.

\proclaim{3. Lemma}
For every $v\in V$ there is $Qv\in V_0$.
\endproclaim
\demo{Proof}
Let $v\in V$ and let $\alpha$ be an 1-form such that $\alpha|V_0=0$.
Then by virtue of Lemma 1 we have
$$
0=(\iota_v\omega)\wedge\omega\wedge\alpha=-\alpha(Qv)\theta,
$$
which shows that $\alpha(Qv)=0$.
\enddemo
\end